\documentclass[11pt, twoside]{article}
\usepackage{mathrsfs}
\usepackage{amsfonts}
\usepackage{graphicx}

\usepackage{amsmath,amssymb}
\usepackage{multicol}
\usepackage{color}
\usepackage{cite}

\newenvironment{proof}[1][ ]{\indent}{\hfill $\Box$}
\newcommand\red{\textcolor[rgb]{0.00,0.00,0.00}}

\allowdisplaybreaks \numberwithin{equation}{section} \makeatletter
\begin{document}

\title{\textbf{An inverse problem of identifying the radiative coefficient in
a degenerate parabolic equation}\thanks{%
Supported by the National Natural Science Foundation of China (Nos.
11061018 and 11261029), Youth Foundation of Lanzhou Jiaotong
University (2011028), Long Yuan Young Creative Talents Support
Program (No. 252003), and the Joint Funds of NSF of Gansu Province
of China (No. 1212RJZA043).}}
\author{Zui-Cha Deng\thanks{%
\textcolor[rgb]{0.00,0.00,0.00}{Corresponding author.}},~~~~Liu Yang \\
Department of Mathematics, Lanzhou Jiaotong University\\
Lanzhou, Gansu, 730070, People's Republic of China \\
E-mail:~~zc\_~deng78@hotmail.com;~~~~l\_~yang218@163.com
   \\
} \maketitle

\begin{abstract}
This work investigates an inverse problem of determining the
radiative coefficient in a degenerate parabolic equation from the
final overspecified data. Being different from other inverse
coefficient problems in which the principle coefficients are assumed
to be strictly positive definite, the mathematical model discussed
in the paper belongs to the second order parabolic equations with
non-negative characteristic form, namely that there exists
degeneracy on the lateral boundaries of the domain. The uniqueness
of the solution is obtained by the contraction mapping principle.
Based on the optimal control framework, the problem is transformed
into an optimization problem and the existence of the minimizer is
established. After the necessary conditions which must be satisfied
by the minimizer are deduced, the uniqueness and stability of the
minimizer are proved. By minor modification of the cost functional
and some \emph{a-priori} regularity conditions imposed on the
forward operator, the convergence of the minimizer for the noisy
input data is obtained in the paper. The results obtained in the
paper are interesting and useful, and can be extended to more
general degenerate parabolic equations.


\end{abstract}

\section{\textbf{\ Introduction }}

\noindent

In this paper, we study an inverse problem of identifying the
radiative coefficient in a degenerate parabolic equation from the
final overspecified data. Problems of this type have important
applications in several fields of applied science and engineering.
The problem can be stated in the following form:

\textbf{Problem P}~~~~Consider the following parabolic equation:
\begin{equation}
\label{1} \left\{\begin{array}{ll} u_t-(a(x)u_x)_x+q(x)u=0,
\hspace{0.3cm} &
(x,t)\in Q=(0,l)\times(0,T],\\
u|_{t=0}=\phi(x), & x\in(0,l),
\end{array}\right.
\end{equation}
where $a$ and $\phi$ are two given smooth functions which satisfy
\begin{equation}\label{p1}
    a(0)=a(l)=0,~~~~a(x)>0,~~x\in(0,l),
\end{equation}
and
\begin{equation}\label{p2}
    \phi(x)\ge 0,~~\phi(x)\not\equiv 0,~~x\in(0,l),
\end{equation}
and $q(x)$ is an unknown coefficient in (\ref{1}). In this paper, we
always assume that $a(x)$ is at least $C^1$ continuous, i.e.,
$a(x)\in C^1[0,l]$. Assume that an additional condition is given as
follows:
\begin{equation}
u(x,T)=g(x),\hspace{0.8cm}x\in[0,l],\label{2}
\end{equation}
where $g$ is a known function. We shall determine the functions $u$
and $q$ satisfying \eqref{1} and \eqref{2}.

If the principle coefficient $a(x)$ is required to be strictly
positive, i.e.,
$$
a(x)\ge a_0>0, ~~x\in[0,l],
$$
then the equation should be rewritten as an initial-boundary value
problem, e.g., the homogeneous Dirichlet boundary value problem as
follows:
\begin{equation}
\label{p3} \left\{\begin{array}{ll} &u_t-(a(x)u_x)_x+q(x)u=0,
\hspace{0.3cm}
(x,t)\in Q,\\
&u|_{x=0}=u|_{x=l}=0, \\
&u(x,0)=\phi(x),
\end{array}\right.
\end{equation}
which is often referred as the classical parabolic equation. The
mathematical model \eqref{p3} arises in various physical and
engineering settings. If \eqref{p3} is used to describe the heat
transfer system, the coefficient $q(x)$ is called the radiative
coefficient which is often dependent on the medium property.

Being different from the ordinary parabolic equation \eqref{p3},
system \eqref{1} belongs to the second order differential equations
with non-negative characteristic form. The main character of such
kinds of equations is degeneracy. It can be easily seen that at
$x=0$ and $x=l$, Eq. \eqref{1} degenerates into two hyperbolic
equations
\begin{eqnarray*}
&& \frac{\partial u}{\partial t}-a'(0)\frac{\partial u}{\partial x}+q(0)u=0, \\
&& \frac{\partial u}{\partial t}-a'(l)\frac{\partial u}{\partial
x}+q(l)u=0.
\end{eqnarray*}
By the well known Fichera's theory (see \cite{R3}) for degenerate
parabolic equations, we know that whether or not boundary conditions
should be given at the degenerate boundaries is determined by the
sign of the Fichera function.

Consider the following second order equation:
\begin{equation}\label{p4}
    \sum_{i,j=1}^m a_{ij}(x)\frac{\partial^2u}{\partial x_i \partial x_j}
    +\sum_{i=1}^m b_i(x)\frac{\partial u}{\partial x_i}+c(x)u=f(x),
\end{equation}
where $x=(x_1,x_2,\cdots,x_m)\in \Omega \subset \mathbb{R}^m$ and
$a_{ij}$ satisfies
$$
a_{ij}=a_{ji},~~i,j=1,2,\cdots,m,
$$
and
$$
\sum_{i,j=1}^m a_{ij}(x)\xi_i \xi_j\ge 0,~~~~\forall x\in
\bar{\Omega},~\xi=(\xi_1,\xi_2,\cdots,\xi_m)\in \mathbb{R}^m.
$$
Let $n=(n_1,n_2,\cdots,n_m)$ be the unit inward normal vector of the
boundary $\partial \Omega$, and let $\partial
\Omega=\Gamma_1\cup\Gamma_2\cup\Gamma_3,$ where $\Gamma_1$ is the
non-characteristic part of $\partial \Omega$, i.e.,
\begin{eqnarray*}
  && \sum_{i,j=1}^m a_{ij}(x)n_i n_j> 0,~~x\in \Gamma_1, \\
  && \sum_{i,j=1}^m a_{ij}(x)n_i n_j= 0,~~x\in \Gamma_2\cup\Gamma_3.
\end{eqnarray*}
Define the following Fichera function:
$$
B(x)=\sum_{i=1}^m \left[b_i(x)-\sum_{j=1}^m\frac{\partial
a_{ij}(x)}{\partial x_j}\right]n_i,
$$
and on $\Gamma_2$ and $\Gamma_3$ it satisfies
$$
B(x)\left \{ \begin{array}{ll} \ge 0, &x\in \Gamma_2, \\
< 0, &x\in \Gamma_3. \end{array} \right .
$$
Then, to guarantee the well-posedness of problem \eqref{p4} one
should give some boundary conditions on $\Gamma_1\cup\Gamma_3$,
while on $\Gamma_2$, they must not be given. For problem \eqref{1},
by denoting $x_1=x$ and $x_2=t$ we have
\begin{eqnarray*}
&& a_{11}=a(x),~~a_{12}=a_{21}=a_{22}=0, \\
&& b_1=a'(x),~~b_2=-1.
\end{eqnarray*}
On the boundary $x=0$, the unit inward normal vector is $(1,0)$. By
direct calculations we have
$$
B(0,t)=b_1-\frac{d}{dx}a_{11}=0.
$$
From the Fichera's theory we know that boundary conditions should
not be given on $x=0$. By analogous arguments, we can also obtain
that on $x=1$ and $t=T$ boundary conditions should not be given,
while on $t=0$ they are indispensable. Therefore, the parabolic
problem \eqref{1} is well defined.

In general, most physical and industrial phenomenons can be
described by the classical parabolic model, such as Eq. \eqref{p3}.
However, with the development of the modern financial mathematics,
more and more degenerate elliptic or parabolic equations arising in
derivatives pricing have to be taken into account. For example, the
well known Black-Scholes equation:
\begin{equation}\label{p7}
    \frac{\partial V}{\partial
t}+\frac{1}{2}\sigma^2(S)S^2\frac{\partial^2 V}{\partial
S^2}+(r-q)S\frac{\partial V}{\partial S}-rV=0,~~~~~~(S,t)\in
[0,\infty)\times[0,T),
\end{equation}
is such the case, where the degenerate parabolic boundary is $S=0$.

For a given coefficient $q(x)$, the degenerate parabolic equation
\eqref{1} which is referred as a direct problem consists of the
determination of the solution from the given initial condition. It
is well known that in all cases the inverse problem is ill-posed or
improperly posed in the sense of Hadamard, while the direct problem
is well-posed (see \cite{b14,b2}). The ill-posedness, particularly
the numerical instability, is the main difficulty for problem
\textbf{P}. Since data errors in the extra condition $g(x)$ are
inevitable, arbitrarily small changes in $g(x)$ may lead to
arbitrarily large changes in $q(x)$, which may make the obtained
results meaningless (see, e.g., \cite{b1,R2}).

Inverse coefficient problems for parabolic equations are well
studied in the literature. However, most of these inverse problems
are governed by classical parabolic equations in which the principle
coefficients are assumed to be strictly positive definite. The
inverse problem of identifying the diffusion coefficient $a(x)$ in
the following parabolic equation
$$
u_t-\nabla\cdot(a(x)\nabla u)=f(x,t),~~~~(x,t)\in \Omega\times(0,T)
$$
from some additional conditions has been investigated by several
authors, e.g., in \cite{101,102,103,104}. In \cite{101,102}, the
output least-squares method with Tikhonov regularization is applied
to the inverse problem and the numerical solution is obtained by the
finite element method. The determination of $a(x)$ with two Neumann
measured data
$$
a(0)u_x(0,t)=k(t),~~~~a(1)u_x(1,t)=h(t),~~t\in[0,T]
$$
has been considered carefully in \cite{103} by the semigroup
approach. In \cite{104}, the inverse problem is reduced to a
nonlinear equation and the uniqueness, as well as the conditional
stability of the solution is proved.

The inverse problem of identification of the radiative coefficient
$q(x)$ in the following heat conduction equation
$$
u_t-\Delta u+q(x)u=0, ~~~~(x,t)\in Q,
$$
from the final overdetermination data $u(x,T)$ has been considered
by several authors, e.g., in \cite{b7,b13,b23,b30,R10}. Moreover,
treatments on the case of purely time dependent $q=q(t)$ can be
found in \cite{k12,k13,b17,b18}. For the general case in which the
unknown coefficient(s) depend(s) on both spatial and temporal
variables, we refer the readers to the references, e.g., in
\cite{R6,b21,MVK,sa}.

Compared with classical parabolic equations, the main difficulty for
degenerate equations lies in the degeneracy of the principle
coefficients which may lead to the corresponding solution has no
sufficient regularity, even if the initial value and the
coefficients are sufficiently smooth functions. Many effective
tools, e.g., the Schauder's type a-priori estimate which has been
extensively applied in classical parabolic equations, are no longer
applicable for the degenerate parabolic equations. The documents
concerned with inverse degenerate problems are quite few in contrast
with those dealt with non-degenerate problems. In \cite{C1}, the
authors investigate an inverse problem of determining the source
term $g$ in the following degenerate parabolic equation
$$
u_t-(x^{\alpha}u_x)_x=g,~~~~(x,t)\in (0,1)\times (0,T),
$$
where $\alpha\in [0,2)$. The uniqueness and Lipschitz stability of
the solution are obtained by the global Carleman estimates which is
introduced in \cite{C2} in 1998. Recently, in \cite{C3} analogous
methods are applied to a nonlinear inverse coefficient problem
arising in the field of climate evolution, where the diffusion
coefficient is assumed to vanish at both extremities of the domain.
For other topics of degenerate parabolic equations, e.g., the null
controllability, we may refer the reader to \cite{C4,C5,C6} and the
reference therein.

The most important inverse problem in which the underlying model is
degenerate may be the reconstruction of local volatility in the
Black-Scholes equation \eqref{p7}. In \cite{b6,b10}, the inverse
problem of identifying the implied volatility $\sigma=\sigma(S)$
from current market prices of options has been considered carefully.
Based on the optimal control framework, the existence, uniqueness of
$\sigma(S)$ and a well-posed algorithm are obtained. Similar results
are derived in \cite{b58}, where a new extra condition, i.e., the
average option premium, is assumed to be known. In \cite{R1}, on the
basis of the parameter-to-solution mapping, the stability and
convergence of approximations for $\sigma(S)$ are gained by Tikhonov
regularization.

It should be mentioned that the degeneracy in the Black-Scholes
equation can be removed by some change of variable (see \cite{R1}).
However, the degeneracy in our problem can not be removed by any
method, which is also the main difficulty in the paper.

To our knowledge, this paper is the first one concerning uniqueness,
stability and convergence of optimal solution in inverse problem for
degenerate parabolic equations such as \eqref{1}.  In this paper, we
use an optimal control framework (see, e.g., \cite{b6,R6,b30,MCS})
to discuss problem \textbf{P} mainly from the theoretical analysis
angle. The outline of the manuscript is as follows: In Section 2,
the uniqueness of the solution for problem \textbf{P} is obtained by
the contraction mapping principle. In Section 3, the inverse problem
\textbf{P} is transformed into an optimal control problem
\textbf{P1} and the existence of minimizer of the cost functional is
proved. The necessary condition of the minimizer is established in
Section 4. By assuming $T$ is relatively small, the local uniqueness
and stability of the minimizer are shown in Section 5. The
convergence of the minimizer with noisy input data is obtained in
Section 6 by some \emph{a-priori} regularity conditions imposed on
the forward operator. Section 7 ends this paper with concluding
remarks.

\section{\textbf{\ Inverse Problem \textbf{P}}}

\noindent

Let's introduce the following function space:
$$
W^{k,\infty}(\Omega)=\{u(x)|~D^{\alpha}u\in
L^{\infty}(\Omega),~~\forall~|\alpha|\leq k\}.
$$

To discuss the uniqueness of the solution, we shall first establish
the weak maximum principle. We would like to consider the more
general equation:
\begin{equation}
\label{CCC6} \left\{\begin{array}{ll} u_t-(a(x)u_x)_x+q(x)u=f(x,t),
\hspace{0.3cm} &
(x,t)\in Q=(0,l)\times(0,T],\\
u|_{t=0}=\phi(x), & x\in(0,l),
\end{array}\right.
\end{equation}

Let
$$G=(0,l)~~~~{\rm and}~~~~G_{\delta}=(-\delta,l+\delta),$$
 where $\delta$ is an arbitrarily small positive constant. Assume that the functions
 $a,~f,~\phi$ and $q$ satisfy the following conditions:
\begin{itemize}
\item $a\in W^{k+1,\infty}(G_{\delta}),~f\in W^{k,\infty}(G_{\delta}\times (0,T)),~\phi\in W^{k+2,\infty}(G_{\delta})$,~$q\in W^{k,\infty}(G_{\delta})$;
\item $a\ge 0,~~~~\forall x\in G_\delta$;
\item $q\ge q_0>0.$
\end{itemize}

\textbf{Theorem 2.1.} (\emph{see \cite{R3}}) Under the above
assumptions, there exists a unique solution $u(x,t)\in
W^{k,\infty}(\bar{Q})$ to the equation \eqref{CCC6}.

For the extra condition \eqref{2} we shall assume that
\begin{equation}
g(x)\in L^{\infty}(0,l). \label{prr}
\end{equation}

\vspace{10pt}

\textbf{Remark 2.1.} The condition $q(x)\ge q_0>0$ is not essential.
In fact, for the case of $q(x)$ with lower bound, i.e., $q(x)\ge
c_0,~c_0<0$, we can make the following function transformation
$$
v(x,t)=u(x,t){\rm e}^{(c_0-1)t}.
$$
One can easily check that $v$ satisfies
$$
v_t-(a(x)v_x)_x+[q(x)-c_0+1]v=\tilde{f}(x,t),
$$
where $q(x)-c_0+1>0$.

\vspace{10pt}

It is known that for classical parabolic equations whose leading
coefficients are assumed to be positive definite, the maximum or
minimum of the solution can only be attained on the parabolic
boundaries, which is known as the famous weak maximum principle.
Such kind of principle is also applicable for degenerate parabolic
equations. The main difference lies in that the maximum or minimum
cannot be attained on the degenerate boundaries.

\vspace{10pt}

Denote
$$
\mathcal{L}u=u_t-(a(x)u_x)_x+q(x)u.
$$

\vspace{10pt}

\textbf{Lemma 2.2.} (\emph{weak maximum principle}) Assuming that
$u(x,t)\in C^2(Q)\cap C(\bar{Q})$, and satisfies $\mathcal{L}u\leq
0$, then we assert that $u(x,t)$ can only attain the positive
maximum at the boundary
$$
\{(x,t)|~~t=0,~x\in [0,l]\}.
$$

\begin{proof}
\textbf{Proof. } Firstly, we assume $u(x,t)$ attains its positive
maximum $M$ at the internal point $P_0(x_0,t_0)\in Q$, i.e.,
$$
u(x_0,t_0)=\max_{\bar{Q}} u(x,t)=M>0.
$$
Then we have
\begin{eqnarray*}
& \displaystyle{\left.\frac{\partial u}{\partial x}\right|_{P_0}=0},       & \hspace{2.5cm}  0<x_0<l,\\
& \displaystyle{\left.\frac{\partial^2 u}{\partial
x^2}\right|_{P_0}\leq 0},
&\hspace{2.5cm} 0<x_0<l, \\
& \displaystyle{\left.\frac{\partial u}{\partial t}\right|_{P_0}=0},  & \hspace{2.5cm} {\rm \hbox{as}}~~t_0<T, \\
& \displaystyle{\left.\frac{\partial u}{\partial t}\right|_{P_0}\ge
0}, & \hspace{2.5cm}{\rm \hbox{as}}~~t_0=T.
\end{eqnarray*}

Hence, we have
$$
\mathcal{L}u|_{P_0}\ge q(x_0)M \ge q_0M>0,
$$
which contradicts with the assumption of the Lemma.

Next, we illustrate that $u(x,t)$ cannot attain its positive maximum
at the degenerate boundaries
$$
\{x=0,~~0<t\leq T\}\bigcup \{x=l,~~0<t\leq T\}.
$$
Without loss of generality, we assume that $u(x,t)$ attains its
positive maximum $M$ at the point $P_1(0,t_1)$, $0<t_1\leq T$. Then
we have£º
\begin{eqnarray*}
& \displaystyle{\left.\frac{\partial u}{\partial t}\right|_{P_1}=0},  & \hspace{2.5cm} {\rm \hbox{µ±}}~~0<t_1<T, \\
& \displaystyle{\left.\frac{\partial u}{\partial t}\right|_{P_1}\ge
0}, & \hspace{2.5cm}{\rm \hbox{µ±}}~~t_1=T.
\end{eqnarray*}
Moreover, from
$$
a(0)=0,~~a(x)>0, ~0<x<l,
$$
we know $a'(0)\ge 0$. Noting that $u(0,t_1)=M$ is a maximum, we have
$$
u(x,t_1)\leq u(0,t_1),~~~~\forall ~x\in (0,l),
$$
which implies
$$
u_x(0,t_1)=\lim_{x\rightarrow 0+}\frac{u(x,t_1)-u(0,t_1)}{x}\leq 0.
$$

Therefore, we have
\begin{eqnarray*}
\mathcal{L}u|_{P_1} &=& \left.\frac{\partial u}{\partial t}\right|_{P_1}-a(0)\left.\frac{\partial^2 u}{\partial x^2}\right|_{P_1}-a'(0)\left.\frac{\partial u}{\partial x}\right|_{P_1}+q(0)u|_{P_1} \\
&\ge & q(0)M\ge q_0M>0,
\end{eqnarray*}
which also contradicts with the assumption of the Lemma.

For same arguments, we know $u(x,t)$ cannot attain its positive
maximum at the degenerate boundary $\{x=l,~~0<t\leq T\}$ either.

This completes the proof of Lemma 2.2.
\end{proof}

\vspace{10pt}

\textbf{Corollary 2.3.} Assuming that $u(x,t)\in C^2(Q)\cap
C(\bar{Q})$ and satisfies $\mathcal{L}u\ge 0$, then $u(x,t)$ can
only attain its negative minimum at the boundary
$$
\{(x,t)|~~t=0,~x\in [0,l]\}.
$$

\vspace{10pt}

\textbf{Theorem 2.4.} Assuming that $u(x,t)\in C^2(Q)\cap
C(\bar{Q})$ is the solution of \eqref{CCC6}, then we have for
$u(x,t)$ the following estimate
\begin{equation}\label{C14}
\max_{\bar{Q}}|u|\leq \max
\left\{\frac{1}{q_0}\sup_{Q}|f|,~~\sup_{[0,l]}|\phi|\right\}.
\end{equation}

\begin{proof}
\textbf{Proof. } Let
$$
M=\left\{\frac{1}{q_0}\sup_{Q}|f|,~~\sup_{[0,l]}|\phi|\right\},
$$
and $v=M\pm u$. One can easily deduce
\begin{eqnarray*}
&&\mathcal{L}v=\mathcal{L}M\pm \mathcal{L}u=Mq(x)\pm f(x,t)\ge 0, \\
&&v|_{t=0}=M\pm \phi(x)\ge 0.
\end{eqnarray*}

From Corollary 2.3, we have
$$
v(x,t)\ge 0,~~~~(x,t)\in \bar{Q}.
$$

Hence
$$
\max_{\bar{Q}} |u|\leq M.
$$

This completes the proof of Theorem 2.4.
\end{proof}

\vspace{10pt}

Now, we will consider the uniqueness of the inverse problem
\textbf{P}.

Let $\mathcal{S}$ be the following function set:
\begin{equation}
\mathcal{S}=\{q(x)|~q(x)\in W^{3,\infty}(G_{\delta}),~q(x)\ge 0\},
\label{Cpr7}
\end{equation}
and assume $a(x),~\phi(x)$ satisfy the following regularity
conditions:
\begin{itemize}
\item $a(x)\in W^{4,\infty}(G_{\delta}),~\phi(x)\in W^{5,\infty}(G_{\delta})$;
\item $a(x)\ge 0,~~~~\forall x\in G_\delta$.
\end{itemize}

Here we make a mini modification for the lower bound of $q(x)$.
Noting $f(x,t)\equiv 0$ and remark 2.1, such the change is
irrelevant.

In the set $\mathcal{S}$, we introduce the following partial order
"$\ge$", i.e., for any $q_1,q_2\in \mathcal{S}$, we call $q_1\ge
q_2$ if and only if
$$
q_1(x)\ge q_2(x),~~~~\forall x \in G_{\delta}.
$$
The partial order "$\leq$" can be defined analogously.

For any $q\in \mathcal{S}$, we introduce a subset of $\mathcal{S}$
denoted by $\mathcal{S}_{q}$ in which all the elements are required
to satisfy the partial order "$\ge$" or "$\leq$" for $q$, i.e.,
$$
\mathcal{S}_{q}=\{q_1\in\mathcal{S}|~~~~q_1\ge q~~\hbox{or}~~q_1\leq
q\}.
$$

\vspace{10pt}

To obtain the uniqueness, we define the following mapping
$\mathbb{P}$ by
\begin{equation}\label{pr8}
    \mathbb{P}[q]=q+\lambda(u(x,T;q)-g(x)),
\end{equation}
where $u(x,t;q)$ is the solution Eq. \eqref{1} with the given
coefficient $q(x)$ and $\lambda>0$ is an adjusting parameter. It can
be easily seen that the existence of fixed points of $\mathbb{P}$ is
equivalent to solutions of the overposed initial value problem.

From Theorem 2.1, we know $u(x,t)\in W^{3,\infty}(\bar{Q})$, which
by the Sobolev embedding theorem implies $u(x,t)\in
C^{2,1}(\bar{Q})$. Then, from Theorem 2.4 we have for $u(x,t)$ the
following estimate:
$$
0\leq u(x,t;q)\leq \|\phi\|_{\infty}.
$$

For any $q,h\in \mathcal{S}$ and $h\ge 0$, one can easily compute
the G\^{a}teaux derivative of $\mathbb{P}$ to obtain
\begin{equation}\label{pr9}
    \mathbb{P}'[q]\cdot h=h-\lambda\hat{u}(x,T;q,h),
\end{equation}
where $\hat{u}(x,t;q,h)$ satisfies the following degenerate
parabolic equation:
\begin{equation}
\label{pr10} \left\{\begin{array}{ll}
&\hat{u}_{t}-(a(x)\hat{u}_x)_x+q(x)\hat{u}=h(x)u,\hspace{0.8cm}
(x,t)\in Q,\\
&\hat{u}|_{t=0}=0. \\
\end{array}\right.
\end{equation}
Noting $h\in \mathcal{S}$ and $u\in W^{3,\infty}(\bar{Q})$, from
Theorem 2.1 we know that there exists a unique weak solution
$\hat{u}\in W^{3,\infty}(\bar{Q})$ to Eq. \eqref{pr10}. For Eq.
\eqref{pr10}, we have from Theorem 2.4 that
$$
0\leq \hat{u}(x,t;q,h)\leq v(x,t),
$$
where the equality on the left if and only if $h\equiv 0$, and
$v(x,t)$ satisfies the following equation:
\begin{equation}
\label{Cpr102} \left\{\begin{array}{ll} &v_{t}-(a(x)v_x)_x+q(x)v=h(x)\|\phi\|_{\infty},\\
&v|_{t=0}=0. \\
\end{array}\right.
\end{equation}

Therefore, for any $\lambda>0$, the righthand side of \eqref{pr9} is
strictly less than $h(x)$. By choosing $\lambda$ sufficiently small
such that $\lambda \|\phi\|_{\infty}<2$, we have
\begin{equation}
\label{Cpr19} \|\mathbb{P}^{'}[q]\cdot
h\|_{\infty}=\|h-\lambda\hat{u}\|_{\infty}<\|h\|_{\infty}.
\end{equation}
From \eqref{Cpr19}, we have for any $q\in \mathcal{S},~ q_1\in
\mathcal{S}_q$,
\begin{equation}
\label{Cpr20}
\|\mathbb{P}[q]-\mathbb{P}[q_1]\|_{\infty}<\|q-q_1\|_{\infty},
\end{equation}
which indicates that if the mapping $\mathbb{P}$ has a fixed point,
then it must be unique in corresponding partial order set.

\vspace{10pt}

\textbf{Theorem 2.5.} If there exist a solution $q(x)\in
\mathcal{S}$ and $u(x,t;q)$ to \eqref{1}/\eqref{2}, then the
solution is unique in the set $\mathcal{S}_q$.

\vspace{10pt}

\textbf{Remark 2.2.} It is well known that the Carleman estimate is
an effective tool to derive uniqueness and conditional stability for
inverse problems (see \cite{C2}). But unfortunately, it fails in
treating the terminal control problems such as inverse problem
\textbf{P}. To the authors' knowledge, the uniqueness obtained in
the paper is so far the best result one can expect.

\section{\textbf{\ Optimal Control Problem }}

\noindent

We have obtained the uniqueness of the solution for problem
\textbf{P} in the previous section. Now we would like to discuss the
regularization of problem \textbf{P}. Before this, we would like to
discuss the forward problem \eqref{CCC6} and give some basic
definitions, lemmas and estimations.

\vspace{10pt}

\textbf{Definition 3.1.} Define $\mathcal{B}$ to be the closure of
$C_0^{\infty}(Q)$ under the following norm:
$$
\|u\|_{\mathcal{B}}^2=\iint_{Q}a(x)(|u|^2+|\nabla u|^2)dxdt,~~~~u\in
\mathcal{B}.
$$

\vspace{10pt}

\textbf{Definition 3.2.} A function $u(x,t)$ is called the weak
solution of \eqref{CCC6}, if $u\in C([0,T];L^2(0,l))\cap
\mathcal{B}$, and for any $\psi\in L^{\infty}((0,T);L^2(0,l))\cap
\mathcal{B}$, $\frac{\partial \psi}{\partial t}\in L^2(Q),$
$\psi(\cdot,T)=0$, the following integration identity holds
\begin{equation}\label{Cc8}
   \iint_{Q}\left(-u\frac{\partial \psi}{\partial t}+a\nabla u\cdot\nabla \psi+qu\psi\right)dxdt-\int_0^l\phi(x)\psi(x,0)dx=\iint_{Q}f\psi dxdt.
\end{equation}

\vspace{10pt}

\textbf{Remark 3.1.} Assuming $u\in C([0,T];L^2(0,l))\cap
\mathcal{B}$, and $\displaystyle{\frac{\partial u}{\partial t}\in
L^2(Q)}$, then \eqref{Cc8} can be rewritten as
$$
 \iint_{Q}\left(\frac{\partial u}{\partial t}\psi+a\nabla u\cdot\nabla \psi+qu\psi\right)dxdt=\iint_{Q}f\psi dxdt,
$$
where $u$ satisfies $\displaystyle{u|_{t=0}=\phi(x)}$ in the sense
of trace.

\vspace{10pt}

\textbf{Theorem 3.1.} For any given $f\in L^\infty(Q),~\phi\in
L^\infty(0,l)$, there exists a unique weak solution
 to \eqref{CCC6} and satisfies the following estimate
$$
\|u\|_{L^{\infty}((0,T),L^2(0,l))}+\|a|\nabla u|^2\|_{L^1(Q)}\leq
C\left(\|f\|^2_{L^2(Q)}+\|\phi\|^2_{L^2(0,l)}\right).
$$
Furthermore, if $a|\nabla \phi|^2\in L^1(0,l)$, then $\frac{\partial
u}{\partial t}\in L^2(Q)$ and
$$
\left\|\frac{\partial u}{\partial t}\right\|_{L^2(Q)}\leq
C\left(\|f\|_{L^2(Q)}+\|\phi\|_{L^2(0,l)}+\|a|\nabla
\phi|^2\|_{L^1(0,l)}\right).
$$

\begin{proof}
\textbf{Proof.} Firstly, we prove the existence. For any given
$0<\varepsilon<1,$ we consider the following regularized problem£º
\begin{equation}
\label{Cc9} \left\{\begin{array}{ll} &\frac{\partial
u_{\varepsilon}}{\partial
t}-(a_{\varepsilon}(x)u_{\varepsilon,x})_x+q(x)u_{\varepsilon}=f(x,t),
\hspace{0.3cm}
(x,t)\in Q,\\
&u_{\varepsilon}(0,t)=u_{\varepsilon}(l,t)=0, \\
&u_{\varepsilon}(x,0)=\phi(x), \\
\end{array}\right.
\end{equation}
where
$$
a_{\varepsilon}(x)=a(x)+\varepsilon,~~x\in [0,l].
$$
From the well-known theory for parabolic equations (see \cite{br8}),
there exists a unique weak solution $u_{\varepsilon}(x,t)$ to Eq.
\eqref{Cc9}.

Then, we will give some \emph{apriori} estimates for
$u_{\varepsilon}(x,t)$. Without loss of generality, we assume that
$u_{\varepsilon}(x,t)$ is the classical solution of \eqref{Cc9}.
Otherwise, one can smooth the coefficients of \eqref{Cc9} and then
consider the solution of the approximation problem.

Multiplying on both sides of \eqref{Cc9} with $u_{\varepsilon}$ and
integrating on $Q_t=[0,l]\times (0,t)$, we have
$$
\iint_{Q_t}\frac{\partial u_{\varepsilon}}{\partial
t}u_{\varepsilon}dxdt-\iint_{Q_t}(a_{\varepsilon}u_{\varepsilon,x})_xu_{\varepsilon}dxdt+\iint_{Q_t}qu_{\varepsilon}^2dxdt=
\iint_{Q_t}fu_{\varepsilon}dxdt.
$$

Integration by parts, we get
\begin{eqnarray}
&& \int_0^l\frac{1}{2}u^2_{\varepsilon}dx+\int_0^t\int_{0}^la_{\varepsilon}|u_{\varepsilon,x}|^2dxdt+\int_0^t\int_{0}^lqu_{\varepsilon}^2dxdt \nonumber \\
&\leq&
\int_0^l\frac{1}{2}\phi^2dx+\frac{1}{2}\int_0^t\int_{0}^l|u_{\varepsilon}|^2dxdt+\frac{1}{2}\int_0^t\int_{0}^lf^2dxdt.
\label{C10}
\end{eqnarray}

From \eqref{C10} and the Gronwall inequality, we have
$$
\max_{0<t\leq
T}\int_0^lu^2_{\varepsilon}dx+\iint_{Q_t}a_{\varepsilon}|u_{\varepsilon,x}|^2dxdt\leq
C\left(\int_0^l\phi^2dx+\iint_{Q_t}f^2dxdt\right).
$$

On the other hand, if $a|\nabla \phi|^2\in L^1(0,l)$, then by
multiplying $\frac{\partial u_{\varepsilon}}{\partial t}$ on both
sides of \eqref{Cc9} and integrating on $Q_t$, we obtain
\begin{eqnarray*}
&&\iint_{Q_t}\left|\frac{\partial u_{\varepsilon}}{\partial t}\right|^2dxdt-\iint_{Q_t}(a_{\varepsilon}u_{\varepsilon,x})_x\cdot\frac{\partial u_{\varepsilon}}{\partial t}dxdt+\iint_{Q_t}qu_{\varepsilon}\frac{\partial u_{\varepsilon}}{\partial t}dxdt \\
&=&\iint_{Q_t}f\frac{\partial u_{\varepsilon}}{\partial t}dxdt.
\end{eqnarray*}
Integrating by parts, we have
\begin{eqnarray}
&&\iint_{Q_t}\left|\frac{\partial u_{\varepsilon}}{\partial
t}\right|^2dxdt
+\iint_{Q_t}\frac{q}{2}\frac{\partial }{\partial t}(u^2_{\varepsilon})dxdt \nonumber \\
&& -\iint_{Q_t}\left[\frac{\partial }{\partial
x}\left(a_{\varepsilon}\frac{\partial u_{\varepsilon}}{\partial
x}\cdot\frac{\partial u_{\varepsilon}}{\partial t}\right)
-a_{\varepsilon}\frac{\partial u_{\varepsilon}}{\partial
x}\cdot\frac{\partial^2 u_{\varepsilon}}{\partial x\partial t}
\right]dxdt \nonumber \\
&=&\iint_{Q_t}\left|\frac{\partial u_{\varepsilon}}{\partial
t}\right|^2dxdt +\iint_{Q_t}\frac{q}{2}\frac{\partial }{\partial
t}(u^2_{\varepsilon})dxdt
+\iint_{Q_t}\frac{a_{\varepsilon}}{2}\frac{\partial }{\partial
t}\left|\frac{\partial u_{\varepsilon}}{\partial
x}\right|^2dxdt \nonumber \\
&=&\iint_{Q_t}f\frac{\partial u_{\varepsilon}}{\partial t}dxdt.
\label{Cc11}
\end{eqnarray}
From \eqref{Cc11}, we get
\begin{eqnarray}
&&\iint_{Q_t}\left|\frac{\partial u_{\varepsilon}}{\partial t}\right|^2dxdt+\int_{0}^l a_{\varepsilon}\left|\frac{\partial u_{\varepsilon}}{\partial x}(\cdot,t)\right|^2dx+\int_{0}^l\frac{q}{2}u^2_{\varepsilon}(\cdot,t)dx  \nonumber \\
&\leq &\int_{0}^l a_{\varepsilon}\phi^2_x dx+\frac{1}{2}\int_{0}^l
q\phi^2dx+\frac{1}{2}\iint_{Q_t}f^2dxdt
+\frac{1}{2}\iint_{Q_t}\left|\frac{\partial
u_{\varepsilon}}{\partial t}\right|^2dxdt. \label{Cc12}
\end{eqnarray}

From \eqref{Cc12}, we have
\begin{equation*}
   \left\|\frac{\partial u_{\varepsilon}}{\partial t}\right\|_{L^2(Q)}\leq C\left(\|f\|_{L^2(Q)}+\|\phi\|_{L^2(0,l)}+\|a_{\varepsilon}|\nabla \phi|^2\|_{L^1(0,l)}\right).
\end{equation*}

Moreover, it follows from the maximum principle that
$$
\|u_{\varepsilon}\|_{L^{\infty}(Q)}\leq C.
$$

From the estimations above, it can be derived that there exists a
subsequence of $\{u_{\varepsilon}\}$ (denoted by itself) and
$$
u\in C([0,T];L^2(0,l)),~~~~\frac{\partial u}{\partial t}\in L^2(Q),
$$
such that
\begin{eqnarray*}
& u_{\varepsilon}\rightarrow u       & \hspace{2.5cm} {\rm in}~~L^2(Q), \\
& \nabla u_{\varepsilon}\rightharpoonup \nabla u & \hspace{2.5cm} {\rm in}~~L_{\rm loc}^2(Q), \\
& \frac{\partial u_{\varepsilon}}{\partial t}\rightharpoonup \frac{\partial u}{\partial t} & \hspace{2.5cm}{\rm in}~~L^2(Q), \\
& a_{\varepsilon}\nabla u_{\varepsilon}\rightharpoonup a\nabla u
&\hspace{2.5cm} {\rm in}~~L^2(Q).
\end{eqnarray*}

Letting $u=u_{\varepsilon}$ in \eqref{Cc8}, we have
\begin{equation*}
   \iint_{Q}\left(-u_{\varepsilon}\frac{\partial \psi}{\partial t}+a\nabla u_{\varepsilon}\cdot\nabla \psi+qu_{\varepsilon}\psi\right)dxdt-\int_0^l\phi(x)\psi(x,0)dx=\iint_{Q}f\psi dxdt.
\end{equation*}

Letting $\varepsilon\rightarrow 0,$ one can immediately obtain
\begin{equation*}
   \iint_{Q}\left(-u\frac{\partial \psi}{\partial t}+a\nabla u\cdot\nabla \psi+qu\psi\right)dxdt-\int_0^l\phi(x)\psi(x,0)dx=\iint_{Q}f\psi
   dxdt,
\end{equation*}
which implies the existence of weak solutions.

Next, we prove the uniqueness of weak solutions. Suppose that
$u_1,~u_2$ be two solutions of \eqref{CCC6} and letÁî
$$
U(x,t)=u_1(x,t)-u_2(x,t),~~~~(x,t)\in Q.
$$
It can be easily seen that $U\in C([0,T];L^2(0,l))\cap \mathcal{B}$,
and for any $\psi\in L^{\infty}((0,T);L^2(0,l))\cap \mathcal{B}$,
$\frac{\partial \psi}{\partial t}\in L^2(Q),$ $\psi(\cdot,T)=0$, the
following integration identity holds
\begin{equation}\label{Cc13}
   \iint_{Q}\left(-U\frac{\partial \psi}{\partial t}+a\nabla U\cdot\nabla \psi+qU\psi\right)dxdt=0.
\end{equation}

For any given $g\in C_0^{\infty}(Q)$, by the existence obtained
above we know that there exists a weak solution $v\in
L^{\infty}((0,T);L^2(0,l))\cap \mathcal{B}$ and $\frac{\partial
v}{\partial t}\in L^2(Q)$ for the following equation
\begin{eqnarray*}
&&-\frac{\partial v}{\partial t}-(a(x)v_x)_x+q(x)v=g(x,t),~~~~(x,t)\in Q, \\
&&v(x,T)=0,~~~~x\in (0,l).
\end{eqnarray*}
Letting $\psi=v$ in \eqref{Cc13}, we obtain
$$
\iint_{Q}Ugdxdt=0.
$$
Noting the arbitrariness of $g$, we have
$$
U(x,t)=0,~~~~{\rm a.e.}~(x,t)\in Q,
$$
i.e.,
$$
u_1(x,t)=u_2(x,t),~~~~{\rm a.e.}~(x,t)\in Q.
$$

This completes the proof of Theorem 3.1.
\end{proof}

\vspace{10pt}

\textbf{Remark 3.2.} The weak solution defined above is on the whole
domain $Q$. If we only consider the spatial case, we can modify the
definition 3.1 as:

\textbf{Definition 3.1'.} Define $\mathscr{H}^1(0,l)$ to be the
closure of $C_0^{\infty}(0,l)$ under the following norm:
$$
\|v\|_{\mathscr{H}^1}^2=\int_{0}^{l}a(x)(|v|^2+|\nabla
v|^2)dxdt,~~~~v\in \mathscr{H}^1(0,l).
$$
For the case of $f\equiv 0$, the definition 3.2 can also be
rewritten as :

\textbf{Definition 3.2'.} A function $u\in H^1((0,T);L^2(0,l))\cap
L^2((0,T);\mathscr{H}^1(0,l)$) is called the weak solution of
\eqref{CCC6}, if $u$ satisfies
\begin{equation}
u(x,0)=\phi(x),~~~~x\in (0,l), \label{Cp9}
\end{equation}
and the following integration identity
\begin{equation}
\int_0^l u_t\psi dx+\int_0^l a\nabla u\cdot \nabla \psi dx+\int_0^l
qu\psi dx=0,~~~~\forall~\psi\in L^2(0,l)\cap \mathscr{H}^1(0,l)
\label{Cp10}
\end{equation}
holds for \emph{a.e.}, $t\in (0,T]$. Then, by analogously arguments,
one can establish the existence, uniqueness and regularity for such
kind of weak solution, which are similar to those of Theorem 3.1.

\vspace{10pt}

\textbf{Remark 3.3.} We recall that the principle coefficient
$a(x)\in C^1[0,1]$. Due to the degeneracy at $x=0$ and $x=l$, from
$u\in \mathscr{H}^1(0,l)$ one can only derive $u\in H^1_{loc}(0,l)$
rather than $u\in H^1(0,l)$, which is different from the case of
non-degenerate. However, we may derive
\begin{equation}\label{t80}
    au_x\rightarrow 0,~~~~{\rm as}~x\rightarrow 0.
\end{equation}
In fact, if \eqref{t80} is not true, i.e., $au_x\rightarrow
k,~k\neq0,$ then we have $u_x\sim \frac{k}{a(x)}$ in
$B_{\delta}(0)\cap [0,l]$, where $B_{\delta}(0)$ is a ball with
$\delta$-radius centered at $x=0$. Note that
\begin{equation*}
    a(x)=a(0)+a'(\xi)x=a'(\xi)x,~~~~\xi\in [0,x],~~x\in B_{\delta}(0)\cap
    [0,l].
\end{equation*}
Hence,
\begin{equation*}
    a|u_x|^2\sim \frac{k^2}{a(x)}\sim \frac{k^2}{a'(\xi)x},
\end{equation*}
which is contradicts with $a|u_x|^2\in L^1(0,l)$. By analogous
arguments, we have $$au_x\rightarrow 0,~~~~{\rm as}~x\rightarrow
l.$$

It should be mentioned that these conclusions are no longer valid
for $a\bar{\in} C^1[0,l]$. For example, let
\begin{equation}\label{t81}
    a(x)=x^{\alpha}(l-x)^{\beta},~~~~0<\alpha,~\beta<1.
\end{equation}
It can be easily seen that $a|u_x|^2\in L^1(0,l)$ cannot guarantee
$au_x\rightarrow 0,$ as $x$ tends to $0$ or $l$. In some references,
e.g., \cite{C4,C5}, the case \eqref{t81} is called the weak
degeneracy and the boundary conditions are indispensable for
corresponding mathematical model, e.g., we shall replace the
equation \eqref{1} by the following initial-boundary value problem:
\begin{equation*}
 \left\{\begin{array}{ll} &u_t-(a(x)u_x)_x+q(x)u=0,
\hspace{0.8cm}
(x,t)\in Q,\\
&u|_{x=0}=u|_{x=l}=0, \\
&u(x,0)=\phi(x).
\end{array}\right.
\end{equation*}

\vspace{20pt}

Since the inverse problem \textbf{P} is ill-posed, i.e., its
solution depends unstably on the data, we turn to consider the
following optimal control problem \textbf{P1}:

Find $\bar{q}(x)\in \mathcal{A}$ such that:
\begin{equation}
J(\bar{q})=\min_{q\in \mathcal{A}}J(q), \label{3}
\end{equation}
\red{where}
\begin{equation}
J(q)=\frac{1}{2}\int_{0}^{l}\left|u(x,T;q)-g(x)\right|^{2}dx+\frac{N}{2}\int_{0}^{l}\left|\nabla
q\right|^{2}dx, \label{4}
\end{equation}
\begin{equation}
\mathcal{A}=\left\{q(x)~|~~0<\alpha \leq q\leq
\beta,~~\|q\|_{H^1(0,l)}<\infty\right\}, \label{5}
\end{equation}
$u(x,t;q)$ is the solution of Eq.(\ref{1}) for a given coefficient
$q(x)\in \mathcal{A}$, $N$ is the regularization parameter, and
$\alpha,~\beta$ are two given positive constants.

\vspace{10pt}

From \eqref{prr} and Theorem 3.1, it can be easily seen that the
control functional \eqref{4} is well-defined for any $q\in
\mathcal{A}$.

\vspace{10pt}

We are now going to show the existence of minimizers to the problem
\eqref{3}. Firstly, we assert that the functional $J(q)$ is of some
continuous property in $\mathcal{A}$ in the following sense.

\vspace{10pt}

\textbf{Lemma 3.2. } For any sequence $\{q_n\}$ in $\mathcal{A}$
which converges to some $q\in \mathcal{A}$ in $L^1(0,l)$ as
$n\rightarrow \infty$, we have
\begin{equation}
\lim_{n\rightarrow \infty}\int_0^l |u(q_n)(x,T)-g(x)|^2dx=\int_0^l
|u(q)(x,T)-g(x)|^2dx. \label{q6}
\end{equation}

\begin{proof}
\textbf{Proof. } \textbf{Step 1:} By taking $q=q_n$ and choosing the
test function as $u(q_n)(\cdot,t)$ in \eqref{Cp10} and then
integrating with respect to $t$, we derive that
\begin{equation}
\|u(q_n;t)\|^2_{L^2(0,l)}+\int_0^t\int_0^l a|\nabla u(q_n;t)|^2dxdt
+\int_0^t\int_0^l q_n|u(q_n;t)|^2dxdt\leq \|\phi\|_{L^2(0,l)}^2
\label{q7}
\end{equation}
for any $t\in (0,T]$.

From \eqref{q7} we know that the sequence $\{u(q_n)\}$ is uniformly
bounded in the space $L^2((0,T);\mathscr{H}^1(0,l))$. So we may
extract a subsequence, still denoted by $\{u(q_n)\}$, such that
\begin{equation}
u(q_n)(x,t)\rightharpoonup u^*(x,t)\in
L^2((0,T);\mathscr{H}^1(0,l)). \label{q8}
\end{equation}

\vspace{10pt}

\textbf{Step 2:} Prove $u^*(x,t)=u(q)(x,t)$.

By taking $q=q_n$ in \eqref{Cp10} and multiplying both sides by a
function $\eta(t)\in C^1[0,T]$ with $\eta(T)=0$, we have
\begin{equation}
\int_0^l u(q_n)_t\psi \eta(t) dx+\int_0^l a\nabla u(q_n)\cdot \nabla
\psi\eta(t) dx+\int_0^l q_nu(q_n)\psi\eta(t) dx=0. \label{q9}
\end{equation}

Then integrating with respect to $t$, we get
\begin{eqnarray}
  -\int_0^l \phi\eta(0)\psi dx &=& -\int_0^T\int_0^l u(q_n)\psi \eta_t(t)dxdt \nonumber\\
  && +\int_0^T\int_0^l \eta(t)a\nabla u(q_n)\cdot \nabla\psi dxdt+
  \int_0^T\int_0^l \eta(t)q(x)u(q_n)\psi dxdt \nonumber\\
  && +\int_0^T\int_0^l \eta(t)(q_n-q)u(q_n)\psi dxdt. \label{q10}
\end{eqnarray}

Letting $n\rightarrow \infty$ in \eqref{q10} and using \eqref{q8},
we obtain
\begin{eqnarray}
  -\int_0^l \phi\eta(0)\psi dx &=& -\int_0^T\int_0^l u^*\psi \eta_t(t)dxdt
  +\int_0^T\int_0^l \eta(t)a\nabla u^*\cdot \nabla\psi dxdt \nonumber\\
  && +\int_0^T\int_0^l \eta(t)q(x)u^*\psi dxdt.  \label{q11}
\end{eqnarray}
By noticing that \eqref{q11} is valid for any $\eta(t)\in C^1[0,T]$
with $\eta(T)=0$, we have
\begin{equation}
\int_0^l u^*_t\psi dx+\int_0^l a\nabla u^*\cdot \nabla \psi
dx+\int_0^l qu^*\psi dx=0,~~~~\forall~\psi\in \mathscr{H}^1(0,l)
\label{q12}
\end{equation}
and $u^*(x,0)=\phi(x)$.

Therefore, $u^*=u(q)$ by the definition of $u(q)$.

\vspace{10pt}

\textbf{Step 3:} Prove $\|u(q_n)(\cdot,T)-g\|_{L^2(0,l)}\rightarrow
\|u(q)(\cdot,T)-g\|_{L^2(0,l)}$ as $n\rightarrow \infty$.

We rewrite \eqref{Cp10} for $q=q_n$ in the form
\begin{eqnarray}
  &&\int_0^l (u(q_n)-g)_t\psi dx+\int_0^l a\nabla (u(q_n)-g)\cdot \nabla \psi
dx+\int_0^l q_n(u(q_n)-g)\psi dx\nonumber\\
&=& -\int_0^l a\nabla g\cdot \nabla \psi dx-\int_0^l q_ng\psi dx.
\label{q13}
\end{eqnarray}

Taking $\psi=u(q_n)-g$ in \eqref{q13} we have
\begin{eqnarray}
  &&\frac{1}{2} \frac{d}{dt}\|u(q_n)-g\|^2_{L^2(0,l)} +\int_0^l a|\nabla (u(q_n)-g)|^2dx
  +\int_0^l q_n|u(q_n)-g|^2dx\nonumber\\
&=& -\int_0^l a\nabla g\cdot \nabla (u(q_n)-g) dx-\int_0^l
q_ng(u(q_n)-g) dx. \label{q14}
\end{eqnarray}

Similar relations hold for $u(q)$, namely
\begin{eqnarray}
  &&\frac{1}{2} \frac{d}{dt}\|u(q)-g\|^2_{L^2(0,l)} +\int_0^l a|\nabla (u(q)-g)|^2dx
  +\int_0^l q|u(q)-g|^2dx\nonumber\\
&=& -\int_0^l a\nabla g\cdot \nabla (u(q)-g) dx-\int_0^l qg(u(q)-g)
dx. \label{q15}
\end{eqnarray}

Subtracting \eqref{q15} from \eqref{q14} we obtain
\begin{eqnarray}
  && \left\{\int_0^l q_n|u(q_n)-g|^2dx-\int_0^l q|u(q)-g|^2dx\right\}+
  \frac{1}{2} \frac{d}{dt}\|u(q_n)-u(q)\|^2_{L^2(0,l)}\nonumber\\
&=&\int_0^l a\nabla g\cdot \nabla
(u(q)-u(q_n))dx+\int_0^lqg(u(q)-u(q_n))dx \nonumber\\
&&+\int_0^l(q-q_n)g(u(q_n)-g)dx+\int_0^l a\nabla(u(q)-u(q_n))\cdot
\nabla(u(q)+u(q_n)-2g)dx \nonumber\\
&&-\int_0^l\frac{d}{dt}[(u(q)-g)(u(q_n)-u(q))]dx. \label{q16}
\end{eqnarray}

Taking $\psi=u(q_n)-u(q)$ in \eqref{Cp10} we have
\begin{eqnarray}
  && \int_0^l u(q)_t(u(q_n)-u(q))dx \nonumber\\
&=&\int_0^l a\nabla u(q)\cdot \nabla (u(q)-u(q_n))dx+\int_0^l
qu(q)(u(q)-u(q_n))dx. \label{q17}
\end{eqnarray}

Similarly, for $(u(q_n)-u(q))_t(u(q)-g)$ we have
\begin{eqnarray}
  && \int_0^l (u(q_n)-u(q))_t(u(q)-g)dx \nonumber\\
&=&\int_0^l a\nabla (u(q_n)-u(q))\cdot \nabla (g-u(q))dx+\int_0^l
q(u(q_n)-u(q))(g-u(q))dx \nonumber\\
&& +\int_0^l(q_n-q)u(q_n)(g-u(q))dx. \label{q18}
\end{eqnarray}

Substituting \eqref{q17} and \eqref{q18} into \eqref{q16} and after
some manipulations, we derive
\begin{eqnarray}
  && \frac{1}{2} \frac{d}{dt}\|u(q_n)-u(q))\|^2_{L^2(0,l)}+\int_0^l a|\nabla (u(q_n)-u(q))|^2dx \nonumber\\
  && +\left\{\int_0^l q_n|u(q_n)-g|^2dx-\int_0^l
  q|u(q)-g|^2dx\right\} \nonumber\\
&=&2\int_0^l q(u(q_n)-u(q))(u(q)-g)dx+\int_0^l (q-q_n)g(u(q_n)-g)dx
\nonumber\\
&&+\int_0^l(q-q_n)u(q_n)(g-u(q))dx:=A_n. \label{q19}
\end{eqnarray}

Then by rewriting the third term on the left side of \eqref{q19}, we
have
\begin{eqnarray}
  && \frac{1}{2} \frac{d}{dt}\|u(q_n)-u(q)\|^2_{L^2(0,l)}+\int_0^l a|\nabla (u(q_n)-u(q)|^2dx
  +\int_0^l q_n|u(q_n)-u(q)|^2dx \nonumber\\
  &=&A_n+\left\{\int_0^l (q-q_n)|u(q)-g|^2dx-2\int_0^l
  q_n(u(q_n)-u(q))(u(q)-g)dx\right\} \nonumber\\
  &:=& A_n+B_n.\label{qq19}
\end{eqnarray}

Integrating over the interval $(0,t)$ for any $t\leq T$, we get
\begin{equation}\label{q20}
    \frac{1}{2} \|u(q_n;t)-u(q;t)\|^2_{L^2(0,l)}\leq
    \int_0^T|A_n+B_n|dt.
\end{equation}

By the convergence of $\{q_n\}$ and the weak convergence of
$\{u(q_n)\}$, one can easily get
\begin{equation}\label{q21}
    \int_0^T|A_n+B_n|dt\rightarrow 0,~~~~~~{\rm as}~~n\rightarrow\infty.
\end{equation}

Combining \eqref{q20} and \eqref{q21} we have
\begin{equation}\label{q22}
    \max_{t\in[0,T]}\|u(q_n;t)-u(q;t)\|_{L^2(0,l)}\rightarrow 0,~~~~~~{\rm as}~~n\rightarrow\infty.
\end{equation}

On the other hand we have from the H\"{o}lder inequality
\begin{eqnarray}
  && \left|\int_0^l |u(q_n)(\cdot,T)-g|^2dx-\int_0^l |u(q)(\cdot,T)-g|^2dx\right| \nonumber\\
  &\leq& \int_0^l |u(q_n)(\cdot,T)-u(q)(\cdot,T)|\cdot|u(q_n)(\cdot,T)+u(q)(\cdot,T)-2g|dx\nonumber\\
  &\leq& \|u(q_n)(\cdot,T)-u(q)(\cdot,T)\|_{L^2(0,l)}\cdot\|u(q_n)(\cdot,T)+u(q)(\cdot,T)-2g\|_{L^2(0,l)}
  .\label{q23}
\end{eqnarray}

From \eqref{prr}, \eqref{q7}, \eqref{q22} and \eqref{q23} we obtain
$$
\lim_{n\rightarrow \infty}\int_0^l |u(q_n)(x,T)-g(x)|^2dx=\int_0^l
|u(q)(x,T)-g(x)|^2dx.
$$

This completes the proof of Lemma 3.2.
\end{proof}

\vspace{10pt}

\textbf{Theorem 3.3. } There exists a minimizer
$\bar{q}\in\mathcal{A}$ of $J(q)$, i.e.
$$
J(\bar{q})=\min_{q\in \mathcal{A}}J(q).
$$

\begin{proof}
\textbf{Proof. } It is obvious that $J(q)$ is non-negative and thus
$J(q)$ has the greatest lower bound $\inf_{q\in \mathcal{A}}J(q)$.
Let $\{q_{n}\}$ be a minimizing sequence, i.e.,
$$\inf_{q\in\mathcal{A}}J(q)\leq J(q_{n})\leq
\inf_{q\in\mathcal{A}}J(q)+\frac{1}{n},~~n=1,2,\cdots.$$

By noticing that $J(q_{n})\leq C$ we deduce
\begin{equation}\label{qq1}
    \left\Vert \nabla q_{n}\right\Vert _{L^{2}(0,l)}\leq C,
\end{equation}
where $C$ is independent of $n$. Noticing the boundedness of
$\{q_{n}\}$ and \eqref{qq1}, we also have
\begin{equation}\label{qq2}
    \|q_{n}\|_{H^{1}(0,l)}\leq C.
\end{equation}
So we can extract a subsequence, still denoted by $\{q_n\}$, such
that
\begin{equation}\label{qq3}
    q_{n}(x)\rightharpoonup \bar{q}(x)\in H^1(0,l)~~{\rm as}~~n\rightarrow\infty.
\end{equation}

By the Sobolev imbedding theorem (see \cite{h1}) we obtain
\begin{equation}\label{qq4}
    \|q_{n}(x)-\bar{q}(x)\|_{L^1(0,l)}\rightarrow 0~~{\rm as}~~n\rightarrow\infty.
\end{equation}
It can be easily seen that $\{q_n(x)\}\in \mathcal{A}$. So we get as
$n\rightarrow \infty$ that
\begin{equation}\label{qq5}
    q_{n}(x)\rightarrow \bar{q}(x)\in \mathcal{A}
\end{equation}
in $L^1(0,l)$.

Moreover, from \eqref{qq3} we have
\begin{equation}\label{qq6}
    \int_0^l|\nabla \bar{q}|^2dx=\lim_{n\rightarrow \infty}\int_0^l
    \nabla q_n\cdot\nabla \bar{q}dx\leq\lim_{n\rightarrow \infty}\sqrt{\int_0^l
    |\nabla q_n|^2dx\cdot\int_0^l|\nabla \bar{q}|^2dx}.
\end{equation}

From Lemma 3.2 and the convergence of $\{q_n\}$, we know that there
exists a subsequence of $\{q_n\}$, still denoted by $\{q_n\}$, such
that
\begin{equation}\label{qq7}
    \lim_{n\rightarrow \infty}\int_0^l |u(q_n)(x,T)-g(x)|^2dx=\int_0^l
|u(\bar{q})(x,T)-g(x)|^2dx.
\end{equation}

From \eqref{qq5}, \eqref{qq6} and \eqref{qq7}, we get

\begin{eqnarray}
  J(\bar{q})&=&\lim_{n\rightarrow \infty}\int_0^l
    |u(q_n)(x,T)-g(x)|^2dx+\int_0^l|\nabla \bar{q}|^2dx \nonumber\\
  &\leq& \lim_{n\rightarrow
    \infty} J(q_n)=\inf_{q\in \mathcal{A}}J(q). \label{qq8}
\end{eqnarray}

Hence, $J(\bar{q})=\displaystyle{\min_{q\in \mathcal{A}}J(q)}$.

This completes the proof of Theorem 3.3. \end{proof}

\section{\textbf{\ Necessary Condition }}

\noindent

\textbf{Theorem 4.1. } Let $q$ be the solution of the optimal
control problem (\ref{3}). Then there exists a triple of functions
$(u,v;q)$ satisfying the following system:
\begin{equation}
\label{6} \left\{\begin{array}{ll} u_t-(au_x)_x+qu=0, \hspace{0.6cm}
&
(x,t)\in Q,\\
u|_{t=0}=\phi(x), & x\in(0,l),
\end{array}\right.
\end{equation}
\begin{equation}
\label{7} \left\{\begin{array}{ll} -v_t-(av_x)_x+qv=0,
\hspace{0.6cm} &
(x,t)\in Q,\\
v|_{t=T}=u(x,T)-g(x), & x\in(0,l),
\end{array}\right.
\end{equation}
and
\begin{equation}
\int_{0}^{T}\int_{0}^{l}uv(q-h)dxdt-N\int_{0}^{l}\nabla q\cdot\nabla
(q-h)dx\ge0  \label{8}
\end{equation}
for any $h\in\mathcal{A}$.

\begin{proof}
\textbf{Proof. } For any $h\in\mathcal{A},~0\leq\delta\leq1$, we
have
$$q_{\delta}\equiv(1-\delta)q+\delta h\in \mathcal{A}.$$

Then
\begin{equation}
J_{\delta}\equiv
J(q_{\delta})=\frac{1}{2}\int_{0}^{l}|u(x,T;q_{\delta})-g(x)|^{2}dx+\frac{N}{2}\int_{0}^{l}|\nabla
q_{\delta}|^{2}dx. \label{9}
\end{equation}

Let $u_{\delta}$ be the solution to the equation (\ref{1}) with
given $q=q_{\delta}$. Since $q$ is an optimal solution, we have
\begin{equation}
\left.\frac{dJ_{\delta}}{d\delta}\right|_{\delta=0}=\int_{0}^{l}
[u(x,T;q)-g(x)]\left.\frac{\partial
u_{\delta}}{\partial\delta}\right|_{\delta=0} dx+N
\int_{0}^{l}\nabla q\cdot\nabla (h-q)dx\geq 0. \label{10}
\end{equation}

\red{Let} $\tilde{u}_{\delta}\equiv\frac{\partial
u_{\delta}}{\partial\delta},$ \red{direct calculations lead} to the
following equation:
\begin{equation}
\label{11} \left\{\begin{array}{rl} &\frac{\partial }{\partial
t}(\tilde{u}_{\delta})-\frac{\partial}{\partial
x}\left(a\frac{\partial\tilde{u}_{\delta}}{\partial x}\right)+q_{\delta}\tilde{u}_{\delta}=(q-h)u_{\delta},\\
&\tilde{u}_{\delta}|_{t=0}=0.\\
\end{array}\right.
\end{equation}

Let $\xi=\tilde {u}_{\delta}|_{\delta=0}$, then $\xi$ satisfies
\begin{equation}
\label{12} \left\{\begin{array}{rl} &\xi_{t}-(a\xi_{x})_x+q\xi=(q-h)u,\\
&\xi|_{t=0}=0.\\
\end{array}\right.
\end{equation}

From (\ref{10}) we have
\begin{equation}
\int_{0}^{l} [u(x,T;q)-g(x)]\xi(x,T)dx+N \int_{0}^{l}\nabla q\cdot
\nabla (h-q)dx\geq 0. \label{13}
\end{equation}

\red{Let} $\mathcal{L}\xi=\xi_{t}-(a\xi_{x})_x+q\xi$, and
\red{suppose} $v$ is the solution of the following problem:
\begin{equation}
\label{14} \left\{\begin{array}{rl} &\mathcal{L}^{\ast}v\equiv -v_{t}-(av_{x})_x+qv=0,\\
&v(x,T)=u(x,T;q)-g(x).\\
\end{array}\right.
\end{equation}
where $\mathcal{L}^{\ast}$ is the adjoint operator of the operator
$\mathcal{L}$.

By the well known Green's formula, we have
\begin{eqnarray}
  && \int_{0}^{T}\int_{0}^{l}(v\mathcal{L}\xi-\xi\mathcal{L}^{\ast}v)dxdt \nonumber\\
  &=& \int_{0}^{T}\int_{0}^{l}(v\xi_t+\xi v_t)dxdt+\int_{0}^{T}\int_{0}^{l}[\xi(av_x)_x-v(a\xi_x)_x]dxdt  \nonumber\\
  &=&\left.\int_{0}^{l}\xi v\right|_{t=0}^{t=T}dx+\int_{0}^{T}\int_{0}^{l}(a\xi
  v_x-av\xi_x)_xdxdt \nonumber\\
  &=& \int_{0}^{l}\xi(x,T)[u(x,T)-g(x)]dx ,\label{15}
\end{eqnarray}
which implies
\begin{equation}
\label{q25} \int_{0}^{T}\int_{0}^{l}v\mathcal{L}\xi
dxdt=\int_{0}^{l}\xi(x,T)[u(x,T)-g(x)]dx.
\end{equation}

Combining (\ref{13}) and (\ref{q25}), one can easily obtain that
$$
\int_{0}^{T}\int_{0}^{l}uv(q-h)dxdt-N\int_{0}^{l}\nabla q\cdot\nabla
(q-h)dx\ge0.
$$

This completes the proof of Theorem 4.1.
\end{proof}

\vspace{10pt}

\section{\textbf{\ Uniqueness and Stability }}

\noindent

The optimal control problem \textbf{P1} is non-convex. So, in
general one may not expect a unique solution. In fact, it is well
known that the optimization technique is a classical tool to yield
"general solution" for inverse problems without unique solution.
However, we find that if the terminal time $T$ is relatively small,
the minimizer of the cost functional can be proved to be local
unique and stable.

\vspace{10pt}

\emph{\textbf{Throughout this paper, if no specific illustration,
$C$ will be denoted the different constants.}}

\vspace{10pt}

\textbf{Lemma 5.1.}~~~~Supposing $u\in \mathscr{H}^1(0,l)$, we have
for any $k\ge 0$,
\begin{eqnarray*}
  && (u-k)^+=\sup (u-k,0)\in  \mathscr{H}^1,\\
  && (u+k)^-=\sup (-(u+k),0)\in  \mathscr{H}^1.
\end{eqnarray*}
Moreover, for \emph{a.e.} $x\in (0,l)$,
$$
\frac{\partial (u-k)^+}{\partial x}\left\{\begin{array}{ll}
                          \frac{\partial u}{\partial x}, ~~~~& {\rm if}~u>k, \\
                          0, & {\rm if}~u\leq k, \\
                        \end{array}\right.
$$
and
$$
\frac{\partial (u+k)^-}{\partial x}\left\{\begin{array}{ll}
                          0, & {\rm if}~u>-k, \\
                          -\frac{\partial u}{\partial x}, ~~~~& {\rm if}~u\leq-k. \\
                        \end{array}\right.
$$

\begin{proof}
\textbf{Proof.} For $u\in \mathscr{H}^1$, we know
$$
\int_0^l a(|u|^2+|\nabla u|^2)dx< +\infty.
$$
Noting $a(x)>0,~x\in(0,l),$ we have for all $\delta>0$,
$$
u\in \mathscr{H}^1(\delta,l-\delta).
$$
By the definition of weak derivative (see \cite{b37}), it can be
easily seen that
$$
(u-k)^+\in \mathscr{H}^1(\delta,l-\delta)
$$
and for {\rm a.e.} $x\in (\delta,l-\delta)$,
$$
\frac{\partial (u-k)^+}{\partial x}\left\{\begin{array}{ll}
                          \frac{\partial u}{\partial x}, ~~~~& {\rm if}~u>k, \\
                          0, & {\rm if}~u\leq k. \\
                        \end{array}\right.
$$
Then we have
$$
\int_{\delta}^{l-\delta}a\left|\left((u-k)^+\right)_x\right|^2dx=\int_{E_{\delta}}a|u_x|^2dx,
$$
where $E_{\delta}=\{~x\in(\delta, l-\delta)~|~u(x)>k~\}$. Since the
quantity $\int_{E_{\delta}}a|u_x|^2dx$ is bounded from above
$\int_0^l a|u_x|^2dx$ which does not depend on $\delta$, by passing
to the limit as $\delta\rightarrow 0$, we get
$$
\int_{0}^{l}a\left|\left((u-k)^+\right)_x\right|^2dx\leq \int_0^l
a|u_x|^2dx<+\infty.
$$
Moreover, the following inequality
$$
\int_{0}^{l}a\left|(u-k)^+\right|^2dx\leq \int_0^l a|u|^2dx<+\infty
$$
is obvious. Hence, $(u-k)^+\in \mathscr{H}^1$. Similar arguments can
be applied to treat the case of $(u+k)^-.$

This completes the proof of Lemma 5.1.
\end{proof}

\vspace{10pt}

Now, we can give a weak maximum principle for the weak solution of
Eq. \eqref{1}.

\vspace{10pt}

\textbf{Lemma 5.2.}~~~~Supposing $\phi\in L^{\infty}(0,l)\cap
\mathscr{H}^1(0,l)$, then we have for $u$ the following estimate:
\begin{equation}\label{t1}
    \|u\|_{\infty}\leq \|\phi\|_{\infty}.
\end{equation}

\begin{proof}
\textbf{Proof.} Let $k=\|\phi\|_{\infty}$. Multiplying the equation
\eqref{1} by $(u-k)^+$, we get from Lemma 5.1
\begin{equation}\label{t2}
    \int_0^l u_t(u-k)^+dx+\int_0^l
    a\left|\left((u-k)^+\right)_x\right|^2dx=-\int_0^l qu(u-k)^+dx.
\end{equation}
Denoting $E=\{~x\in(0, l)~|~u(x)>k~\}$, one has
\begin{equation}\label{t3}
    -\int_0^l qu(u-k)^+dx=-\int_E qu(u-k)^+dx\leq 0.
\end{equation}
From \eqref{t2} and \eqref{t3}, we have for all $t\in [0,T]$,
\begin{equation*}
    \frac{1}{2}\frac{d}{dt}\int_0^l
    \left|(u-k)^+\right|^2dx=\int_0^l u_t(u-k)^+dx\leq 0,
\end{equation*}
which implies $t\mapsto \|(u-k)^+(t)\|^2_{L^2}$ is decreasing on
$[0,T]$. Since $(\phi-k)^+\equiv 0$, we deduce that for all $t\in
[0,T]$ and for \emph{a.e.} $x\in (0,l),~u(x,t)\leq k$.

By analogous arguments for $(u+k)^-$, we can obtain that for all
$t\in [0,T]$ and for \emph{a.e.} $x\in (0,l),~u(x,t)\ge -k$.

This completes the proof of Lemma 5.2.
\end{proof}

\vspace{10pt}

\textbf{Lemma 5.3.}~~~~For Eq. (\ref{7}) we have the following
estimate:
\begin{equation}
\|v\|_{\infty}\leq \|u(x,T)-g(x)\|_{\infty}. \label{q30}
\end{equation}

\begin{proof}
\textbf{Proof. } Let $\tau=T-t$, then \eqref{7} is reduced to
\begin{equation*}
 \left\{\begin{array}{rl} &v_{\tau}-(av_{
x})_x+qv=0, ~~(x,t)\in Q, \\
&v|_{\tau=0}=u(x,T)-g(x). \\
\end{array}\right.
\end{equation*}

The rest of the proof is similar to that of Lemma 5.2. \end{proof}

\vspace{10pt} Suppose $g_1(x)$ and $g_2(x)$ be two given functions
which satisfy the condition \eqref{prr}. Let $q_{1}(x)$ and
$q_{2}(x)$ be the minimizers of problem \textbf{P1} corresponding to
$g=g_i,~(i=1,2)$ respectively, and $\{u_{i},v_{i}\},(i=1,2)$ be
solutions of system (\ref{6})/(\ref{7}) in which $q=q_{i},(i=1,2)$
respectively.

Setting
$$u_{1}-u_{2}=U,~~v_{1}-v_{2}=V,~~q_{1}-q_{2}=\mathcal{Q},$$
then $U$ and $V$ satisfy
\begin{equation}
\label{a1} \left\{\begin{array}{rl} &U_{t}-(aU_{
x})_x+q_{1}U=-\mathcal{Q}u_2,\\
&U|_{t=0}=0,
\end{array}\right.
\end{equation}
\begin{equation}
\label{a2} \left\{\begin{array}{rl} &-V_{t}-(aV_{
x})_x+q_{1}V=-\mathcal{Q}v_2,\\
&V|_{t=T}=U(x,T)-(g_1-g_2).
\end{array}\right.
\end{equation}

\vspace{10pt}

\textbf{Lemma 5.4. } For any bounded continuous function $k(x)\in
C(0,l)$, we have
$$
\|k\|_{\infty}\leq |k(x_0)|+\sqrt{l}\|\nabla k\|_{L^2(0,l)},
$$
where $x_{0}$ is a fixed point in $(0, l)$.

\begin{proof}
\textbf{Proof. } For $0<x<l$ we have
\begin{eqnarray*}
\left\vert k(x)\right\vert&\leq&\left\vert
k(x_{0})\right\vert+\left\vert \int_{x_{0}}^{x}k^{'}dx\right\vert \\
&\leq&\left\vert
k(x_{0})\right\vert+\left(\int_{0}^{l}dx\right)^{\frac{1}{2}}
\left(\int_{0}^{l}\left\vert \nabla
k\right\vert^{2}dx\right)^{\frac{1}{2}}.
\end{eqnarray*}

This completes the proof of the Lemma 5.4.
\end{proof}

\vspace{10pt} \textbf{Lemma 5.5.}~~~~For equation (\ref{a1}) we have
the following estimate:
\begin{equation}
\max_{0\leq t\leq T}\int_{0}^{l}U^2dxdt \leq C(\max
|\mathcal{Q}|)^{2}\int_{0}^{T}\int_{0}^{l}|u_{2}|^{2}dxdt,
\label{a3}
\end{equation}
where $C$ is independent of $T$.

\begin{proof}
\textbf{Proof. } From equation (\ref{a1}) we have for $0<t\leq T$
\begin{eqnarray}
  && \int_{0}^{l}\int_{0}^{t}\left(\frac{U^2}{2}\right)_{t}dxdt-\int_{0}^{t}\int_{0}^{l}
(aU_{x})_xUdxdt+\int_{0}^{t}\int_{0}^{l}q_{1}U^2dxdt \nonumber\\
  &=& -\int_{0}^{t}\int_{0}^{l}u_2\mathcal{Q}Udxdt. \label{d7}
\end{eqnarray}

Integrating by parts we obtain
\begin{eqnarray}
  && \int_{0}^{l}\left.\frac{U^2}{2}\right|_{(x,t)}dx+\int_{0}^{t}\int_{0}^{l}
aU_{x}^2dxdt-\left.\int_{0}^{t}aU_xU\right|_{x=0}^{x=l}dt+\int_{0}^{t}\int_{0}^{l}q_{1}U^2dxdt \nonumber\\
  &\leq& \int_{0}^{t}\int_{0}^{l}U^2dxdt+(\max|Q|)^2\int_{0}^{t}\int_{0}^{l}|u_2|^2dxdt, \label{d7}
\end{eqnarray}
which implies
\begin{eqnarray}
  && \int_{0}^{l}\left.\frac{U^2}{2}\right|_{(x,t)}dx+\int_{0}^{t}\int_{0}^{l}
aU_{x}^2dxdt \nonumber\\
  &\leq& \int_{0}^{t}\int_{0}^{l}U^2dxdt+(\max|Q|)^2\int_{0}^{t}\int_{0}^{l}|u_2|^2dxdt. \label{q31}
\end{eqnarray}

From the Gronwall inequality and \eqref{q31} we have
$$
\int_{0}^{l}U^2dxdt+\int_{0}^{T}\int_{0}^{l}aU_x^2dxdt \leq C(\max
|\mathcal{Q}|)^{2}\int_{0}^{T}\int_{0}^{l}|u_{2}|^{2}dxdt.
$$

This completes the proof of Lemma 5.5.
\end{proof}

\vspace{10pt} \textbf{Lemma 5.6.}~~~~For equation (\ref{a2}) we have
the following estimate:
\begin{eqnarray}
&& \max_{0\leq t\leq
T}\int_{0}^{l}V^2dx+\int_{0}^{T}\int_{0}^{l}a|V_{x}|^2dxdt \nonumber\\
  &\leq &
  C(\max|\mathcal{Q}|)^2\int_0^T\int_0^l(|u_2|^2+|v_2|^2)dxdt+C\int_{0}^{l}|g_1-g_2|^2dx.
  \label{d11}
\end{eqnarray}
where $C$ is independent of $T$.

\begin{proof}
\textbf{Proof. } From equation (\ref{a2}) we have
\begin{eqnarray*}
  && \int_{t}^{T}\int_{0}^{l}-\left(\frac{V^2}{2}\right)_{t}dxdt-\int_{t}^{T}\int_{0}^{l}(aV_{x})_xVdxdt
+\int_{t}^{T}\int_{0}^{l}q_{1}V^2dxdt \\
  &=&
  -\int_{t}^{T}\int_{0}^{l}v_2\mathcal{Q}Vdxdt.
\end{eqnarray*}

Integrating by parts we obtain that
\begin{eqnarray}
&&\int_{0}^{l}\left.\frac{V^2}{2}\right|_{(x,t)}dx+\int_{t}^{T}\int_{0}^{l}a|V_{x}|^2dxdt
+\int_{t}^{T}\int_{0}^{l}q_{1}V^2dxdt \nonumber\\
&\leq& \int_0^l|U(x,T)|^2dx+\int_{0}^{l}|g_1-g_2|^2dx-\int_{t}^{T}\int_{0}^{l}v_2\mathcal{Q}Vdxdt \nonumber\\
&\leq&
\int_0^l|U(x,T)|^2dx+\int_{0}^{l}|g_1-g_2|^2dx+\int_t^T\int_0^l\frac{V^2}{2}dxdt
\nonumber\\
&& +\frac{1}{2}(\max|\mathcal{Q}|)^2\int_t^T\int_0^l|v_2|^2dxdt.
\label{d12}
\end{eqnarray}

From Lemma 5.5 and \eqref{d12} we have
\begin{eqnarray}
  && \int_{0}^{l}\left.\frac{V^2}{2}\right|_{(x,t)}dx+\int_{t}^{T}\int_{0}^{l}a|V_{x}|^2dxdt \nonumber\\
  &\leq & \int_t^T\int_0^l\frac{V^2}{2}dxdt+\int_{0}^{l}|g_1-g_2|^2dx  \nonumber\\
  &&+C(\max|\mathcal{Q}|)^2\int_0^T\int_0^l(|u_2|^2+|v_2|^2)dxdt.
  \label{d14}
\end{eqnarray}

From the Gronwall inequality we have
\begin{eqnarray*}
&& \max_{0\leq t\leq
T}\int_{0}^{l}V^2dx+\int_{0}^{T}\int_{0}^{l}a|V_{x}|^2dxdt \\
  &\leq &
  C(\max|\mathcal{Q}|)^2\int_0^T\int_0^l(|u_2|^2+|v_2|^2)dxdt+C\int_{0}^{l}|g_1-g_2|^2dx.
\end{eqnarray*}

This completes the proof of Lemma 5.6.
\end{proof}

\vspace{10pt}

\textbf{Theorem 5.7. } Let $q_{1}(x),~q_{2}(x)$ be the minimizers of
the optimal control problem \textbf{P1} corresponding to
$g_{1}(x),g_{2}(x)$, respectively. If there exists a point $x_{0}\in
(0,l)$ such that
$$
q_{1}(x_{0})=q_{2}(x_{0}),
$$
then for relatively small $T$ we have
$$
\max_{x\in(0,l)}|q_1-q_2|\leq \frac{C
l^{\frac{1}{3}}}{N^{\frac{1}{3}}}\|g_1-g_2\|_{L^2(0,l)},
$$
where the constant $C$ is independent of $T$, $l$ and $N$.

\begin{proof}
\textbf{Proof. } \red{By taking} $h=q_{2}$ when $q=q_{1}$ and
\red{taking} $h=q_{1}$ when $q=q_{2}$ in (\ref{8}), we have
\begin{equation}
\int_{0}^{T}\int_{0}^{l}(q_{1}-q_{2})u_{1}v_{1}dxdt-N\int_{0}^{l}\nabla
q_{1}\red{\cdot}\nabla (q_{1}-q_{2})dx\ge0,  \label{a5}
\end{equation}
\begin{equation}
\int_{0}^{T}\int_{0}^{l}(q_{2}-q_{1})u_{2}v_{2}dxdt-N\int_{0}^{l}\nabla
q_{2}\red{\cdot}\nabla (q_{2}-q_{1})dx\ge0,  \label{a6}
\end{equation}
where $\{u_{i},v_{i}\},(i=1,2)$ are solutions of system
(\ref{6})/(\ref{7}) with $q=q_{i}~(i=1,2)$, respectively.

From (\ref{a5}) and (\ref{a6}) we have
\begin{eqnarray}
 && N\int_{0}^{l}|\nabla (q_{1}-q_{2})|^2dx \nonumber\\
 &\leq&
\int_{0}^{T}\int_{0}^{l}(u_{1}v_{1}-u_{2}v_{2})(q_{1}-q_{2})dxdt \nonumber\\
&\leq&
\int_{0}^{T}\int_{0}^{l}(u_{1}v_{1}-u_{2}v_{1}+u_{2}v_{1}-u_{2}v_{2})(q_{1}-q_{2})dxdt
\nonumber\\
&\leq&
\int_{0}^{T}\int_{0}^{l}\mathcal{Q}v_1Udxdt+\int_{0}^{T}\int_{0}^{l}\mathcal{Q}u_2Vdxdt.
\label{a7}
\end{eqnarray}

From the assumption of Theorem 5.7, there exists a point $x_{0}\in
(0,l)$ such that
\begin{equation}
\mathcal{Q}(x_{0})=q_{1}(x_{0})-q_{2}(x_{0})=0. \label{a16}
\end{equation}

From Lemma 5.4 and \eqref{a16} we have
\begin{equation}
\max_{x\in(0,l)}\left\vert \mathcal{Q}(x)\right\vert\leq
\sqrt{l}\left(\int_{0}^{l}\left\vert \nabla
\mathcal{Q}\right\vert^{2}dx\right)^{\frac{1}{2}}. \label{a8}
\end{equation}

From (\ref{a7}), (\ref{a8}) and the Young inequality , we obtain
that
\begin{eqnarray}
\max|\mathcal{Q}|^2 &\leq& l\int_{0}^{l}\left\vert \nabla
\mathcal{Q}\right\vert^{2}dx \nonumber \\
&\leq& \frac{l}{N}\int_{0}^{T}\int_{0}^{l}\mathcal{Q}(Uv_1+Vu_2)dxdt
\nonumber \\
&\leq & \frac{1}{2l}\int_{0}^{l}\left\vert
\mathcal{Q}\right\vert^{2}dx+\frac{Tl^2}{2N^2}\int_{0}^{T}\int_{0}^{l}|Uv_1+Vu_2|^2dxdt
\nonumber \\
& \leq&
\frac{1}{2}\max|\mathcal{Q}|^2+\frac{Tl^2}{N^2}\|v_1\|^2_{\infty}\int_{0}^{T}\int_{0}^{l}U^2dxdt
\nonumber \\
&& +\frac{Tl^2}{N^2}\|u_2\|^2_{\infty}\int_{0}^{T}\int_{0}^{l}V^2dxdt \nonumber\\
 & \leq&
\frac{1}{2}\max|\mathcal{Q}|^2+C\frac{T^2l^2}{N^2}\|v_1\|^2_{\infty}\cdot\left(\int_{0}^{T}\int_{0}^{l}|u_2|^2dxdt\right)\cdot\max|\mathcal{Q}|^2 \nonumber \\
&&+C\frac{T^2l^2}{N^2}\|u_2\|^2_{\infty}\cdot\left(\int_{0}^{T}\int_{0}^{l}(|u_{2}|^{2}+|v_2|^2)dxdt\right)\cdot\max|\mathcal{Q}|^2 \nonumber \\
&& +C\frac{T^2l^2}{N^2}\int_{0}^{l}|g_1-g_2|^2dx,\label{a9}
\end{eqnarray}
where we have used estimates \eqref{a3} and \eqref{d11}.

From Lemmas 5.2 and 5.3 we have
\begin{eqnarray}
\|v_1\|_{\infty},~\|v_2\|_{\infty},~\|u_2\|_{\infty}\leq C.
\label{d22}
\end{eqnarray}

From (\ref{a9}) and (\ref{d22}) we have
\begin{equation}
\max|\mathcal{Q}|^2\leq
C\frac{T^3l^2}{N^2}\max|\mathcal{Q}|^2+C\frac{T^2l^2}{N^2}\int_{0}^{l}|g_1-g_2|^2dx.
\label{a11}
\end{equation}

\red{Choose} $T<<1$ such that
\begin{equation}
C\frac{T^3l^2}{N^2}=\frac{1}{2}. \label{a12}
\end{equation}

Combining (\ref{a11}) and (\ref{a12}) one can easily get
\begin{equation}
\max_{x\in(0,l)}|q_1-q_2|\leq \frac{C
l^{\frac{1}{3}}}{N^{\frac{1}{3}}}\|g_1-g_2\|_{L^2(0,l)}. \label{a13}
\end{equation}

This completes the proof of the Theorem 5.7. \end{proof}

\vspace{20pt}

\textbf{Remark 5.1.}  It should be mentioned that the regularization
parameter plays a major role in the numerical simulation of
ill-posed problems. From Theorem 5.7 we can obtain that if there
exists a constant $\delta$ such that
$$
 \|g_1-g_2\|\leq
 \delta,~~~~\text{and}~~~~\frac{\delta^2}{N^{\frac{2}{3}}}\rightarrow 0,
 $$
 then the reconstructed optimal solution is unique and stable, which is consistent with the existed results
 (see, for instance, \cite{R2}). Note that the estimate \eqref{a13}
 is based on \eqref{a12} from which we can see
 $T=O(N^{\frac{2}{3}})$. Since the parameter $N$ is often taken to
 be very small, particularly in numerical computations, Theorem 5.7
 is indeed the local well-posedness of the optimal solution. For more detailed
  discussion on the regularization parameter, we refer the readers to the
  references, e.g., in \cite{R2,KK3}.

\section{Convergence Analysis}

\noindent

In this section, we would like to discuss the convergence of the
optimal solution. It has been shown in previous section that the
optimal solution is stable and unique, which is very important in
numerical process. However, the optimization problem is just a
"modified problem" rather than the original one. Therefore, it is
necessary to investigate what about the difference between the
optimal solution of the optimization problem and the exact solution
of the original problem.

We assume that the "real solution" $g(x)$ is attainable, i.e., there
exists a $q^*(x)\in H^1(0,l)$ such that
\begin{equation}\label{t11}
    u(x,T;q^*)=g(x),
\end{equation}
and that an upper bound $\delta$ for the noisy level
\begin{equation}\label{t12}
    \|g^{\delta}-g\|_{L^2(0,l)}\leq \delta,
\end{equation}
of the observation is known \emph{a priori}.

It should be mentioned that for terminal control problems, it is
rather difficult to derive the convergence. To the authors'
knowledge, there is no any convergence results for the optimal
control problem with the cost functional whose form is similar to
\eqref{4}.

In this paper, we introduce the following auxiliary control problems
with observations averaged over the given terminal time interval
$[T-\sigma,T]$:
\begin{equation}\label{t13}
    J_{\sigma}(q)=\frac{1}{2\sigma}\int_{T-\sigma}^{T}\int_{0}^{l}\left|u(x,t;q)-g(x)\right|^{2}dxdt+\frac{N}{2}\int_{0}^{l}\left|\nabla
q\right|^{2}dx.
\end{equation}
Note that as $\sigma\rightarrow 0^+$,
\begin{equation*}
    \frac{1}{2\sigma}\int_{T-\sigma}^{T}\int_{0}^{l}\left|u(x,t;q)-g(x)\right|^{2}dxdt\rightarrow
    \int_{0}^{l}\frac{1}{2}\left|u(x,T;q)-g(x)\right|^{2}dx,
\end{equation*}
which implies $J_{\sigma}(q)\rightarrow J(q)$. Analogously, instead
of \eqref{t12}, we assume that for the real solution $q^*(x)$,
\begin{equation}\label{t14}
    \frac{1}{2\sigma}\int_{T-\sigma}^{T}\int_{0}^{l}\left|u(x,t;q^*)-g^{\delta}(x)\right|^{2}dxdt\leq
    \frac{1}{2}\delta^2.
\end{equation}

Define the following forward operator $u(q)$:
\begin{eqnarray}
  u(q):~~~~~~~~~~~~~~~~~\mathcal{A} &\rightarrow& H^1((0,T);L^2(0,l))\cap L^2((0,T);\mathscr{H}^1(0,l)) \nonumber\\
  u(q)(x,t) &=&u(x,t;q(x)), \nonumber
\end{eqnarray}
where $u(x,t;q(x))$ is the solution of the variational problem
\eqref{Cp10} for $q\in \mathcal{A}$. For any $q\in \mathcal{A}$ and
$p\in H^1(0,l)$, one can easily deduce that the G\^{a}teaux
directional differential $u'(q)p$ satisfies a homogeneous initial
condition and solves
\begin{equation}\label{t16}
    \int_0^l \left(u'(q)p\right)_t\varphi dx+\int_0^l
    a\nabla\left(u'(q)p\right)\cdot\nabla\varphi dx+\int_0^l qu'(q)p\varphi
    dx=-\int_0^l pu(q)\varphi dx,
\end{equation}
for any $\varphi\in L^2(0,l)\cap\mathscr{H}^1(0,l)$. For the
remainder term $R(q)=u(p+q)-u(q)-u'(q)p$, we have the following
variational characterization.

\vspace{10pt}

\textbf{Lemma 6.1.} For any $q\in\mathcal{A}$ and $p\in H^1(0,l)$
such that $p+q\in \mathcal{A}$, the remainder
$R(q)=u(p+q)-u(q)-u'(q)p$ solves
\begin{equation}\label{t17}
    \int_0^l \left(R(q)\right)_t\varphi dx+\int_0^l
    a\nabla\left(R(q)\right)\cdot\nabla\varphi dx+\int_0^l qR(q)\varphi
    dx=\int_0^l p\varphi\left(u(q)-u(q+p)\right) dx,
\end{equation}
for any $\varphi\in L^2(0,l)\cap \mathscr{H}^1(0,l)$.

\begin{proof}
\textbf{Proof.} Note that $u(q+p)$ satisfies
\begin{equation}\label{t18}
    \int_0^l \left(u(q+p)\right)_t\varphi dx+\int_0^l
    a\nabla\left(u(q+p)\right)\cdot\nabla\varphi dx+\int_0^l (q+p)u(q+p)\varphi
    dx=0.
\end{equation}
Subtracting \eqref{t18} from \eqref{Cp10} and denoting
$W=u(q+p)-u(q)$, we obtain
\begin{equation}\label{t19}
    \int_0^l\varphi W_t dx+\int_0^l
    a\nabla W\cdot\nabla\varphi dx+\int_0^l qW\varphi
    dx=-\int_0^l pu(q+p)\varphi dx.
\end{equation}
Now \eqref{t17} follows by subtracting \eqref{t16} from \eqref{t19}.

This completes the proof of Lemma 6.1.
\end{proof}

\vspace{10pt}

To obtain the convergence,  we shall require some source conditions.
We introduce the following linear operator $F(q)$:
\begin{eqnarray}
  F(q):~~~~~L^2((0,T);L^2(0,l)) &\rightarrow& L^2(0,l) \nonumber\\
  F(q)\Phi &=&-\frac{1}{\sigma}\int_{T-\sigma}^T u(q)\Phi dt, ~~~~\forall \Phi\in L^2((0,T);L^2(0,l)),
  \label{t20}
\end{eqnarray}
where $u(q)$ is the solution of \eqref{Cp10}. Using the equation
\eqref{t16}, we immediately see that for any $p\in H^1(0,l)$ and any
$\varphi\in L^2(0,l)\cap \mathscr{H}^1(0,l)$, the following holds:
\begin{eqnarray}
  <F(q)\varphi,~p> &=& -\frac{1}{\sigma}\int_{T-\sigma}^T pu(q)\varphi dt \nonumber\\
  &=& \frac{1}{\sigma}\int_{T-\sigma}^T\left[\left(u'(q)p\right)_t\varphi +
    a\nabla\left(u'(q)p\right)\cdot\nabla\varphi +qu'(q)p\varphi
    \right]dxdt, \label{t21}
\end{eqnarray}
where $<\cdot~,~\cdot>$ denote the scalar product in $L^2(0,l)$.
Since $\nabla$ is a linear operator, we can define its adjoint
operator $\nabla^*$ by
\begin{equation}\label{t22}
    <\nabla^* \omega,~\varphi>_{L^2(0,l)}=<\omega,~\nabla
    \varphi>_{L^2(0,l)},~~~~\forall~\omega\in H^1(0,l),~\varphi \in
    H^1(0,l).
\end{equation}
It can be easily seen that if $\varphi\in H_0^1(0,l)$, then
$\nabla^*$ is equivalent to $\nabla$. In this paper, we will only
need a weak form of $\nabla^*\nabla$.

\vspace{10pt}

\textbf{Theorem 6.2.} Assume that there exists a function
\begin{equation*}
    \varphi\in H_0^1((T-\sigma,T);L^2(0,l))\cap L^2((T-\sigma,T);\mathscr{H}^1(0,l))
\end{equation*}
such that the following source condition holds in the weak sense:
\begin{equation}\label{t23}
    F(q^*)\varphi=\nabla^*\nabla q^*
\end{equation}
with $F(q^*)$ defined by \eqref{t20}, i.e., for any $p\in H^1(0,l)$,
\begin{equation}\label{t24}
    <F(q^*)\varphi,~p>=<\nabla^*\nabla q^*,~p>=<\nabla q^*,~\nabla
    p>.
\end{equation}
Furthermore, assume that
\begin{equation}\label{t41}
    \nabla\cdot(a\nabla \varphi)\in L^2((T-\sigma,T);L^2(0,l))
\end{equation}
and $q_N^{\delta}$ satisfies
\begin{equation}\label{t42}
    q_N^{\delta}(0)=q^*(0),~~~~~~q_N^{\delta}(l)=q^*(l).
\end{equation}
Then, with $N\sim\delta$, we have
\begin{equation}\label{t25}
    \int_0^l \left|q_N^{\delta}-q^*\right|^2dx\leq C\delta,
\end{equation}
and
\begin{equation}\label{t26}
   \frac{1}{\sigma}\int_{T-\sigma}^{T}\int_{0}^{l}\left|u(q_N^{\delta})-u(q^*)\right|^{2}dxdt\leq
    C\delta^2,
\end{equation}
where $q_N^{\delta}$ is a minimizer of \eqref{t13} with $g$ replaced
by $g^{\delta}$, $u(q_N^{\delta})$ is the solution of the
variational problem \eqref{Cp10} with $q=q_N^{\delta}$, and $C$ is a
positive constant independent of $\delta,~N$ and $T$.

\begin{proof}
\textbf{Proof.} Noting that $q_N^{\delta}$ is a minimizer of
\eqref{t13}, we have
\begin{equation*}
    J_{\sigma}(q_N^{\delta})\leq J_{\sigma}(q^*),
\end{equation*}
which implies
\begin{equation}\label{t27}
    \frac{1}{2\sigma}\int_{T-\sigma}^{T}\int_{0}^{l}\left|u(q_N^{\delta})-g^{\delta}\right|^{2}dxdt+\frac{N}{2}\int_{0}^{l}\left|\nabla
q_N^{\delta}\right|^{2}dx\leq
\frac{1}{2}\delta^2+\frac{N}{2}\int_{0}^{l}\left|\nabla
q^*\right|^{2}dx.
\end{equation}
From \eqref{t27}, we can derive
\begin{eqnarray}
  && \frac{1}{2\sigma}\int_{T-\sigma}^{T}\int_{0}^{l}\left|u(q_N^{\delta})-g^{\delta}\right|^{2}dxdt+\frac{N}{2}\int_{0}^{l}\left|\nabla
q_N^{\delta}-\nabla q^*\right|^{2}dx \nonumber\\
  &\leq & \frac{1}{2}\delta^2+\frac{N}{2}\int_{0}^{l}\left|\nabla
q^*\right|^{2}dx-\frac{N}{2}\int_{0}^{l}\left|\nabla
q_N^{\delta}\right|^{2}dx+\frac{N}{2}\int_{0}^{l}\left|\nabla
q_N^{\delta}-\nabla q^*\right|^{2}dx \nonumber\\
  &=& \frac{1}{2}\delta^2+N\int_{0}^{l}\nabla
q^*\cdot \nabla(q^*-q_N^{\delta})dx \nonumber\\
  &=& \frac{1}{2}\delta^2+N\left<\nabla q^*,~\nabla \left(q^*-q_N^{\delta}\right)\right>. \label{t28}
\end{eqnarray}

Using \eqref{t21} and \eqref{t24}, we have for the last term in
\eqref{t28} that
\begin{eqnarray}
  && \left<\nabla q^*,~\nabla \left(q^*-q_N^{\delta}\right)\right>  \nonumber\\
  &=& \left<F(q^*)\varphi,~q^*-q_N^{\delta}\right> \nonumber\\
  &=& -\frac{1}{\sigma}\int_{T-\sigma}^{T}\int_{0}^{l}(q^*-q_N^{\delta})u(q^*)\varphi dxdt \nonumber\\
  &=& \frac{1}{\sigma}\int_{T-\sigma}^T\int_0^l\left[\left(u'(q^*)(q^*-q_N^{\delta})\right)_t\varphi +
    a\nabla\left(u'(q^*)(q^*-q_N^{\delta})\right)\cdot\nabla\varphi\right.
    \nonumber\\
  && \hspace{40pt}\left.+q^*u'(q^*)(q^*-q_N^{\delta})\varphi
    \right]dxdt \label{t29}
\end{eqnarray}

Let
\begin{equation}\label{t30}
    R_N^{\delta}:=u(q_N^{\delta})-u(q^*)-u'(q^*)(q_N^{\delta}-q^*).
\end{equation}
Using this notation, we obtain
\begin{eqnarray}
  && N\left<\nabla q^*,~\nabla \left(q^*-q_N^{\delta}\right)\right>  \nonumber\\
  &=& \frac{N}{\sigma}\int_{T-\sigma}^T\int_0^l\left[\left(R_N^{\delta}\right)_t\varphi +
    a\nabla R_N^{\delta}\cdot\nabla\varphi+q^*R_N^{\delta}\varphi
    \right]dxdt \nonumber\\
  && -\frac{N}{\sigma}\int_{T-\sigma}^T\int_0^l\left[u(q_N^{\delta})-u(q^*)\right]_t\varphi
    dxdt-\frac{N}{\sigma}\int_{T-\sigma}^T\int_0^la\nabla\left(u(q_N^{\delta})-u(q^*)\right)\cdot\nabla\varphi
    dxdt
    \nonumber\\
  &&-\frac{N}{\sigma}\int_{T-\sigma}^T\int_0^l q^*\left[u(q_N^{\delta})-u(q^*)\right]\varphi
    dxdt \nonumber\\
  &=&I_1+I_2+I_3+I_4.  \label{t31}
\end{eqnarray}

Now, we need to estimate $I_1-I_4$. The main idea is to control
$I_1-I_4$ by the left-side item of inequality \eqref{t28}.

For $I_1$, we use \eqref{t17} to get
\begin{equation}\label{t32}
    I_1=\frac{N}{\sigma}\int_{T-\sigma}^T\int_0^l \left(q_N^{\delta}-q^*\right)\left[u(q^*)-u(q_N^{\delta})\right]\varphi
    dxdt.
\end{equation}
From \eqref{t32} and the H\"{o}lder inequality, we have
\begin{eqnarray}
 |I_1| &\leq& \frac{N}{\sigma}\int_{T-\sigma}^T\int_0^l
 \left|\left(q_N^{\delta}-q^*\right)\varphi\right|\cdot\left|u(q^*)-g^{\delta}\right|
    dxdt  \nonumber\\
  && +\frac{N}{\sigma}\int_{T-\sigma}^T\int_0^l
 \left|\left(q_N^{\delta}-q^*\right)\varphi\right|\cdot\left|g^{\delta}-u(q_N^{\delta})\right|
    dxdt \nonumber\\
  &\leq& \frac{N}{\sigma}\int_{T-\sigma}^T
 \left\|\left(q_N^{\delta}-q^*\right)\varphi\right\|_{L^2(0,l)}\cdot \left\|u(q^*)-g^{\delta}\right\|_{L^2(0,l)}dt
    \nonumber\\
  &&+\frac{N}{\sigma}\int_{T-\sigma}^T
 \left\|\left(q_N^{\delta}-q^*\right)\varphi\right\|_{L^2(0,l)}\cdot \left\|g^{\delta}-u(q_N^{\delta})\right\|_{L^2(0,l)}dt.  \label{t33}
\end{eqnarray}
Using \eqref{5} and the Young inequality, we obtain
\begin{eqnarray}
 |I_1| &\leq& \frac{1}{8\sigma}\int_{T-\sigma}^T\int_0^l
 \left|u(q^*)-g^{\delta}\right|^2dxdt+CN^2\int_{T-\sigma}^T\int_0^l
 \left|\left(q_N^{\delta}-q^*\right)\varphi\right|^2dxdt \nonumber\\
  && +\frac{1}{16\sigma}\int_{T-\sigma}^T\int_0^l
\left|g^{\delta}-u(q_N^{\delta})\right|^2
    dxdt+CN^2\int_{T-\sigma}^T\int_0^l
 \left|\left(q_N^{\delta}-q^*\right)\varphi\right|^2
    dxdt \nonumber\\
  &\leq& \frac{1}{8}\delta^2+\frac{1}{16\sigma}\int_{T-\sigma}^T\int_0^l
\left|g^{\delta}-u(q_N^{\delta})\right|^2
    dxdt+CN^2\int_{T-\sigma}^T\int_0^l
 \left|\varphi\right|^2
    dxdt,  \label{t34}
\end{eqnarray}
where we have used the assumption \eqref{t14}.

For $I_2$, using integration by parts with respect to $t$ and
noticing $\varphi\in H_0^1((T-\sigma,T);L^2(0,l))$, we derive
\begin{eqnarray}
 |I_2|&=&
 \frac{N}{\sigma}\left|\int_{T-\sigma}^T\int_0^l\left(u(q_N^{\delta})-u(q^*)\right)\varphi_t
    dxdt\right| \nonumber\\
 &\leq&
 \frac{N}{\sigma}\int_{T-\sigma}^T\int_0^l\left|\left(u(q_N^{\delta})-u(q^*)\right)\varphi_t\right|
    dxdt \nonumber\\
 &\leq& \frac{N}{\sigma}\int_{T-\sigma}^T\int_0^l\left|\left(u(q_N^{\delta})-g^{\delta}\right)\varphi_t\right|
    dxdt+\frac{N}{\sigma}\int_{T-\sigma}^T\int_0^l\left|\left(g^{\delta}-u(q^*)\right)\varphi_t\right|
    dxdt \nonumber\\
  &\leq& \frac{1}{8}\delta^2+\frac{1}{16\sigma}\int_{T-\sigma}^T\int_0^l
\left|g^{\delta}-u(q_N^{\delta})\right|^2
    dxdt+CN^2\int_{T-\sigma}^T\int_0^l
 \left|\varphi_t\right|^2
    dxdt,  \label{t35}
\end{eqnarray}

For $I_3$, using integration by parts with respect to $x$ and
noticing $a(0)=a(l)=0$, we obtain
\begin{eqnarray}
 |I_3| &=& \frac{N}{\sigma}\left|\int_{T-\sigma}^T\int_0^l a\nabla\left(u(q_N^{\delta})-u(q^*)\right)\cdot\nabla\varphi
    dxdt\right| \nonumber\\
  &=&
  \frac{N}{\sigma}\left|\int_{T-\sigma}^T\left\{\left.a(x)\left(u(q_N^{\delta})-u(q^*)\right)\frac{d\varphi}{dx}\right|_{x=0}^{x=l}
    \right.\right. \nonumber\\
  &&\hspace{60pt}\left.\left.-\int_0^l \left(u(q_N^{\delta})-u(q^*)\right)\nabla\cdot(a\nabla\varphi)dx\right\}dt\right| \nonumber\\
  &\leq& \frac{N}{\sigma}\int_{T-\sigma}^T\int_0^l
  \left|u(q_N^{\delta})-u(q^*)\right|\cdot\left|\nabla\cdot(a\nabla\varphi)\right|dxdt \nonumber\\
  &\leq& \frac{N}{\sigma}\int_{T-\sigma}^T\int_0^l
  \left|u(q_N^{\delta})-g^{\delta}\right|\cdot\left|\nabla\cdot(a\nabla\varphi)\right|dxdt \nonumber\\
  && +\frac{N}{\sigma}\int_{T-\sigma}^T\int_0^l
  \left|g^{\delta}-u(q^*)\right|\cdot\left|\nabla\cdot(a\nabla\varphi)\right|dxdt \nonumber\\
  &\leq& \frac{1}{8}\delta^2+ \frac{1}{16\sigma}\int_{T-\sigma}^T\int_0^l
  \left|u(q_N^{\delta})-g^{\delta}\right|^2dxdt \nonumber\\
  && +CN^2\int_{T-\sigma}^T\int_0^l
 \left|\nabla\cdot(a\nabla\varphi)\right|^2dxdt.  \label{t36}
\end{eqnarray}

The last term $I_4$ can be estimated similarly using the Young
inequality:
\begin{eqnarray}
 |I_4| &\leq& \frac{N}{\sigma}\int_{T-\sigma}^T\int_0^l
 |q^*|\left|u(q_N^{\delta})-g^{\delta}\right||\varphi|
    dxdt \nonumber\\
  &&+\frac{N}{\sigma}\int_{T-\sigma}^T\int_0^l
 |q^*|\left|g^{\delta}-u(q^{*})\right||\varphi|
    dxdt \nonumber\\
  &\leq& \frac{1}{8}\delta^2+ \frac{1}{16\sigma}\int_{T-\sigma}^T\int_0^l
  \left|u(q_N^{\delta})-g^{\delta}\right|^2dxdt
  +CN^2\int_{T-\sigma}^T\int_0^l
 \left|\varphi\right|^2dxdt,  \label{t37}
\end{eqnarray}
where we have used the bound of $q^*$.

Combining \eqref{t28}, \eqref{t31} and \eqref{t34}-\eqref{t37}, we
obtain
\begin{eqnarray}
  && \frac{1}{2\sigma}\int_{T-\sigma}^{T}\int_{0}^{l}\left|u(q_N^{\delta})-g^{\delta}\right|^{2}dxdt+\frac{N}{2}\int_{0}^{l}\left|\nabla
q_N^{\delta}-\nabla q^*\right|^{2}dx \nonumber\\
  &\leq& \frac{1}{2}\delta^2+N\left<\nabla q^*,~\nabla \left(q^*-q_N^{\delta}\right)\right> \nonumber\\
  &\leq& \frac{1}{2}\delta^2+\sum_{j=1}^4 |I_j| \nonumber\\
 &\leq&
 \delta^2+\frac{1}{4\sigma}\int_{T-\sigma}^{T}\int_{0}^{l}\left|u(q_N^{\delta})-g^{\delta}\right|^{2}dxdt \nonumber\\
&& +CN^2\int_{T-\sigma}^T\int_0^l\left(
 \left|\varphi\right|^2+|\varphi_t|^2+\left|\nabla\cdot(a\nabla\varphi)\right|^2\right)dxdt. \label{t38}
\end{eqnarray}
From \eqref{t38} and noticing the regularity of $\varphi$, we have
\begin{equation}\label{t39}
    \frac{1}{4\sigma}\int_{T-\sigma}^{T}\int_{0}^{l}\left|u(q_N^{\delta})-g^{\delta}\right|^{2}dxdt+\frac{N}{2}\int_{0}^{l}\left|\nabla
q_N^{\delta}-\nabla q^*\right|^{2}dx\leq \delta^2+CN^2.
\end{equation}
By choosing $N\sim \delta$, one can easily get
\begin{equation}\label{t40}
    \frac{1}{\sigma}\int_{T-\sigma}^{T}\int_{0}^{l}\left|u(q_N^{\delta})-g^{\delta}\right|^{2}dxdt+N\int_{0}^{l}\left|\nabla
q_N^{\delta}-\nabla q^*\right|^{2}dx\leq C\delta^2.
\end{equation}
The estimate \eqref{t25} follows immediately from \eqref{t40} and
the Poincar\`{e} inequality.

This completes the proof of Theorem 6.2.
\end{proof}

\vspace{10pt}

\textbf{Remark 6.1.} The motivation of replacing the cost functional
\eqref{4} by \eqref{t13} mainly lies in the difficulty in treating
the second integration term in \eqref{t31}. In fact, if we choose
the functional form \eqref{4}, then we can deduce the second term in
\eqref{t31} (denoted by $\tilde{I_2}$) to be
\begin{equation*}
    \tilde{I_2}=-\frac{N}{\sigma}\int_0^l
    \left(u(q_N^{\delta})-u(q^*)\right)_t(\cdot,T)\varphi dx.
\end{equation*}
Since we have no any information regarding to the $t$-derivative of
the real and approximate solution, it is quite difficult, even
impossible, to control the term $\tilde{I_2}$ by the left-hand side
of \eqref{t28}, and thus we cannot obtain any convergence.

\section{Concluding Remarks}

\noindent

The inverse problem of identifying the coefficient in parabolic
equations from some extra conditions is very important in some
engineering texts and many industrial applications. Classical
parabolic models are plentifully discussed and developed well, while
documents dealt with degenerate parabolic models are quite few.

In this paper, we solve the inverse problem \textbf{P} of recovering
the radiative coefficient $q(x)$ in the following degenerate
parabolic equation
$$
u_t-(au_{x})_x+q(x)u=0
$$
in an optimal control framework. Being different from other works
(for example \cite{101,104}) which also treat with inverse radiative
coefficient problems, the mathematical model discussed in the paper
contains degeneracy on the lateral boundaries. Furthermore, unlike
the well known Black-Scholes equation whose degeneracy can be
removed by some change of variable, the degeneracy in our problem
can not be removed by any method. On the basis of the optimal
control framework, the existence, uniqueness, stability and
convergence of the minimizer for the cost functional are
established.

The paper focuses on the theoretical analysis of the 1-D inverse
problem. For the multidimensional case, i.e., the determination of
$q(x)$ in the following equation
$$
u_t-\nabla\cdot(a(x)\nabla u)+q(x)u=0,~~~~(x,t)\in
Q=\Omega\times(0,T],
$$
where the principle coefficient $a(x)$ satisfies
$$
a(x)\ge 0,~~x\in \bar{\Omega}
$$
and $\Omega\subset \mathbb{R}^m~(m\ge1)$ is a given bounded domain,
the method proposed in the paper is also applicable.

\section*{Acknowledgements}
\noindent

This work is supported by the National Natural Science Foundation of
China (Nos. 11061018 and 11261029), Youth Foundation of Lanzhou
Jiaotong University (2011028), Long Yuan Young Creative Talents
Support Program (No. 252003), and the Joint Funds of NSF of Gansu
Province of China (No. 1212RJZA043).

\end{document}